\documentclass{siamart190516}

\usepackage{amsmath,amssymb,enumerate,bm,url}
\usepackage{hyperref,cleveref}
\usepackage{microtype}
\usepackage{arydshln}
\usepackage{appendix}
\usepackage{graphicx}
\usepackage{caption}
\newsiamremark{remark}{Remark}
\usepackage[T1]{fontenc}
\usepackage{lmodern}

\makeatletter

\@ifundefined{showcaptionsetup}{}{ \PassOptionsToPackage{caption=false}{subfig}}
\makeatother

\begin{document}

\title{DeepMartNet - A Martingale Based Deep Neural Network Learning Method for Dirichlet
BVPs and Eigenvalue Problems of Elliptic PDEs in $R^{\MakeLowercase{d}}$ \thanks{December 20, 2023, first version appeared in arxiv preprint arXiv:2311.09456. 2023 Nov 15.}}
\author{ Wei Cai\thanks{Corresponding author, Department of Mathematics, Southern
Methodist University, Dallas, TX 75275 (\texttt{cai@smu.edu})}
\and Andrew He\thanks{Department of Mathematics, Southern
Methodist University, Dallas, TX 75275.}
\and Daniel Margolis\thanks{Department of Mathematics, Southern Methodist
	University, Dallas, TX 75275.}
}

\maketitle

\begin{abstract}
In this paper, we propose DeepMartNet - a Martingale based deep neural network learning method for solving
Dirichlet boundary value problems (BVPs) and eigenvalue problems for elliptic partial differential equations (PDEs) in high
dimensions or domains with complex geometries. The method is based on Varadhan's
Martingale problem formulation for the
BVPs/eigenvalue problems where a loss function enforcing the Martingale property for the PDE solution is used for an
efficient optimization by sampling the stochastic processes associated with corresponding
elliptic operators. High dimensional numerical results for BVPs of the linear and nonlinear Poisson-Boltzmann equation and eigenvalue problems
of the Laplace equation and a Fokker-Planck equation demonstrate the capability of the proposed DeepMartNet learning method in solving high dimensional PDE problems.
\end{abstract}



\section{Introduction}

Computing eigenvalues and/or eigenfunctions for elliptic operators or solving
boundary value problem of partial differential equations (PDEs) in high dimensions and optimal stochastic control problems are among the
key tasks for many scientific computing applications, e.g., ground states and band
structure calculation in quantum systems, complex biochemical systems, communications and manufacturing productions, etc.
Deep neural network (DNN) has been utilized to solve high dimensional PDEs and stochastic control problems due to its capability of approximating high dimensional functions. The first attempt to use DNN to solve high dimensional quasi-linear parabolic PDE was carried out in \cite{DeepBSDE} using Pardoux-Peng theory \cite{peng} of nonlinear Feynman-Kac formulas, which connect the solution of the parabolic PDEs with that of backward stochastic differential equations (SDEs) \cite{peng20}. The loss function for the training of the DeepBSDE uses the terminal condition of the PDEs. The DeepBSDE framework was also applied to solve stochastic control problems \cite{hanEcontrol}. The DeepBSDE started a new approach of SDE based DNN approximations to high dimensional PDEs and stochastic control problems. The work in \cite{raissi18} \cite{zhang22} extended this idea to use a loss function based on the pathwise comparison of two  stochastic processes, one from the BSDE in the Pardoux-Peng theory and one from the PDEs solution, this variant of SDE based DNNs has the potential to find the solution in the whole domain compared with one point solution in the original DeepBSDE.
Recently, a diffusion Monte Carlo DNN method was developed \cite{han20}
using the connection between stochastic processes and the solution of elliptic
equations and the backward Kolmogorov equation to build a loss function for
eigen-value calculations. The capability of SDE paths exploring high dimensional spaces make them a good candidate to be used for DNN learning. It should be mentioned that other approaches to solving high dimensional PDEs include the Feynmann-Kac formula based Picard iteration \cite{Jentzen20} \cite{Epicard} and stochastic dimension gradient descent method \cite{sdgd}.

In addition to the Pardoux-Peng BSDE approach by linking stochastic processes and the solution of PDEs, another powerful probabilistic method is the Varadhan's Martingale problem approach  \cite{Varadhan71, shreve}, which were used to derive a probabilistic weak form for the PDE's solution with a specific Martingale  related to the PDE solution through the Ito formulas \cite{klebaner}. The equivalence of the classic weak solution of boundary value problems (BVPs) of elliptic PDEs using bilinear forms and the probabilistic one were established for the Schrodinger equation for the Neumann BVP \cite{hsu84,hsu85}, and then, for the Robin BVP \cite{Pap1990}.  This Martingale weak form for the PDEs solution and in fact also for a wide class of stochastic control problem \cite{davis05, cohen15} provides a new venue to tackle high dimensional PDE and stochastic control problems \cite{deepMart}, and this paper focuses on the case of high dimensional PDEs only. The Martingale formulation of the PDE solution is a result of Ito formula and the Martingale nature of Ito integrals. As a simple conclusion from the broader Martingale property, the  Feynman-Kac formula provides one point solution of the PDE using expectation of the underlying SDE paths originating from that point. The Martingale based DNN to be studied in this paper, termed DeepMartNet \cite{deepMart}, is trained based on a loss function enforcing the conditional expectation definition of the Martingale. It will be shown with extensive numerical tests that using the same set of SDE paths originating from one single point, the DeepMartNet can in fact provide approximation to solutions of the BVPs and eigenvalue problems of elliptic PDEs over the whole solution doamin in high dimensions. In this sense, the DeepMartNet is able to extract more information for the PDE solutions from the SDE paths originating from one point than the (pre-machine learning) traditional use of the one point solution Feynman-Kac formula.

The rest of this paper is organized as follows. In section 2, a brief review of existing SDE-based DNN methods for solving PDEs is given and Section 3 will present the Martingale problem formulation for the BVP problem of a general  elliptic PDE for the case of third kind Robin boundary condition, which includes both Dirichlet and Neumann BVPs as special limiting cases. Section 4 will present the Martingale based DeepMartNet for solving high dimensional PDE problems.  Numerical results for the Dirichlet BVPs and eigenvalue problems will be present in Section 5. The implementation of the DeepMartNet for the Neumann and Robin boundary conditions will be addressed in a follow-up paper, which will involve reflecting diffusion processes in finite domains and the computation of local times of the processes. Section 6 will present conclusions and future work.



\section{A review of SDE based DNNs for solving PDEs}
To set the background for the Martingale based DNNs, we will first briefly review some existing DNN based on diffusion paths from SDEs.

Let us first consider a terminal value problem for quasi-linear elliptic PDEs
\begin{equation}\label{eq-paraPDE}
\partial_t u + \mathcal{L} u  = \phi,  \quad \textbf{x} \in R^d,
\end{equation}
with a terminal condition $u(T, \mathbf{x}) = g(\mathbf{x})$, where the
differential operator $\mathcal{L}$ is given as
\begin{equation}
\mathcal{L}=\mu^{\top}\nabla+\frac{1}{2}Tr(\sigma\sigma^{\top}\nabla \nabla^{\top}%
) = \mu^{\top}\nabla+ \frac{1}{2}Tr(A \nabla \nabla^{\top}%
), \label{gen}%
\end{equation}
and  $\mu=\mu(t, \textbf{x}, u, \nabla u), \sigma=\sigma(t, \textbf{x}, u, \nabla u) $ and the diffusion coefficient matrix
\begin{equation}
    A= (a_{ij})_{d\times d}= \sigma\sigma^{\top}.
\end{equation}

The aim is to find the solution at $\mathbf{x}, t=0,$ $u(0, \mathbf{x})$, and the solution of \eqref{eq-paraPDE} is related to a coupled FBSDE \cite{peng}
\begin{align}
\begin{split} \label{eq-Xt}
d\mathbf{X}_t &= \mu(t, \textbf{X}_t, Y_t, \textbf{Z}_t) dt + \sigma(t, \textbf{X}_t, Y_t) d\textbf{B}_t,  \\
\textbf{X}_0 &= \xi,
\end{split} \\
\begin{split} \label{eq-Yt}
dY_t &= \phi(t, \textbf{X}_t, Y_t, \textbf{Z}_t) dt + \textbf{Z}^T_t \sigma(t, \textbf{X}_t, Y_t) d\textbf{B}_t, \\
Y_T &= g(\textbf{X}_T),
\end{split}
\end{align}
where $\textbf{X}_t$, $Y_t$ and $\textbf{Z}_t$ are $d$, 1 and $d$-dimensional stochastic processes, respectively, that are adapted to $\{ \mathcal{F}_t: 0 \le t \le T \}$ - the natural filtration from the d-dimensional Brownian motion $\textbf{B}_t$ .
Specifically, we have the following relations,
\begin{equation}
Y_t = u(t, \textbf{X}_t), \quad \textbf{Z}_t = \nabla u(t, \textbf{X}_t).
\end{equation}

\medskip
\noindent {\bf DeepBSDE.} As the first work of using SDEs to train DNN, the Deep BSDE \cite{DeepBSDE}  trains the network with input $\textbf{X}_0 = \textbf{x}$ and output $\textbf{Y}_0=u(0, \mathbf{x})$.
Applying the { Euler--Maruyama scheme} (EM) to the FBSDE (\ref{eq-Xt}) and (\ref{eq-Yt}), respectively, we have
\begin{align}
\textbf{X}_{n+1} & \approx \textbf{X}_n + \mu(t_n, \textbf{X}_n, \textbf{Y}_n, \textbf{Z}_n) \Delta t_n + \sigma(t_n, \textbf{X}_n, \textbf{Y}_n) \Delta \textbf{B}_n, \label{eq-Xn} \\
Y_{n+1} & \approx Y_n + \phi(t_n, \textbf{X}_n, Y_n, \textbf{Z}_n) \Delta t_n + Z^T_n \sigma(t_n, \textbf{X}_n, Y_n) \Delta \textbf{B}_n. \label{eq-Yn}
\end{align}

The missing $\textbf{Z}_{n+1}$ at $t_{n+1}$ will be then approximated by a neural network (NN) with parameters $\theta_n$
\begin{equation}
 \nabla u(t_{n}, \textbf{X}_{n} | \theta_{n}) \sim  \textbf{Z}_{n}= \nabla u(t_{n}, \textbf{X}_{n}).
\end{equation}
Loss function: with an ensemble average approximation, the loss function is defined as
\begin{equation}
Loss_{bsde}(Y_0, \theta)=E \left\lVert u(T, \textbf{X}_T) - g(\textbf{X}_T) \right\rVert^2,
\end{equation}
where
\[
u(T, \textbf{X}_T)=Y_N.
\]

Trainable parameters are $\{ Y_0, \theta_n, n=1, \cdots, N \}$.

\medskip
\noindent {\bf FBSNN.} A forward backward neural network (FBSNNs) proposed in \cite{raissi18} uses the mismatch between two stochastic processes to build the loss function for the DNN, which aims to train a DNN $u_{\theta}(x,t)$ in the whole domain. The following is an improved version of the approach in \cite{zhang22}.

\begin{itemize}

\item Markov chain one. Starting with $\textbf{X}_0 = \textbf{x}$ , $Y_0=u_{\theta}(\textbf{x},0)$,
\begin{align}
\begin{split}
\textbf{X}_{n+1} &= \textbf{X}_n + \mu(t_n, \textbf{X}_n, Y_n, Z_n) \Delta t_n + \sigma(t_n, \textbf{X}_n, Y_n) \Delta \textbf{B}_n, \\
Y_{n+1} & = Y_n + \phi(t_n, \textbf{X}_n, Y_n, \textbf{Z}_n) \Delta t_n + \textbf{Z}^T_n \sigma(t_n, \textbf{X}_n, Y_n) \Delta \textbf{B}_n, \\
\textbf{Z}_{n+1} &{}={} \nabla u(t_{n+1}, \textbf{X}_{n+1}).
\end{split}
\end{align}
\item Markov chain two
\begin{equation}
{ Y^\star_{n+1} = u(t_{n+1}, \textbf{X}_{n+1}).}
\end{equation}
\item The loss function is a Monte Carlo approximation of
\begin{equation}
E \left[ \frac{1}{N} \sum_{n=1}^{N} \left \lVert Y_{n} - Y^\star_{n} \right \rVert^2 + 0.02 \left \lVert Y^\star_{N} - g(\textbf{X}_N) \right \rVert^2 + 0.02 \left \lVert \textbf{Z}_{N} - \nabla g(\textbf{X}_N) \right \rVert^2 \right].
\end{equation}
\end{itemize}
Numerical results of half-order convergence of $u_{\theta}(\textbf{x},t)$, similar to the order of the underlying Euler-Maruyama scheme, is observed.

\medskip
\noindent {\bf Diffusion Monte Carlo DNN eigensolver.}
In another approach similar to the power iteration method, a diffusion-Monte Carlo method was proposed \cite{han20} through a fixed point of semi-group formulation for the eigenvalue problem of the following linear elliptic operator (i.e. $\mu=\mu(\textbf{x}), \sigma=\sigma(\textbf{x})$),
\begin{equation}
    \mathcal{L}\Psi =\lambda \Psi.
    \label{evp0}
\end{equation}
Equation \eqref{evp0} can be reformulated as a virtual time dependent backward parabolic PDE with the sought-after eigenfunction as the terminal condition, i.e.,
\begin{equation}
    \partial_t u(t,\textbf{x}) +\mathcal{L}u(t,\textbf{x})-\lambda u(t,\textbf{x})=0, \quad u(T, \textbf{x})=\Psi(\textbf{x}).
\end{equation}
Thus, the following fix-point property holds,
\begin{equation}
    u(T-t, \cdot)=P^\lambda_t \Psi, \quad P^\lambda_T \Psi=\Psi,
\end{equation}
where the semi-group for the evolutionary system is formally defined by
\begin{equation}
    P^\lambda_t=e^{-(T-t) \mathcal{L}}.
\end{equation}

The discretized backward in time evolution (with time step $\Delta t$ of the backward parabolic equation) mimics the power iteration of the semi-group operator $e^{-n\Delta t \mathcal{L}} $, which will converge to the lowest eigenfunction for $\mathcal{L}$ as in a power method. The loss function is set to be $||P^\lambda_T \Psi-\Psi||^2$ , while the evolution of PDE solution is done by two SDEs, instead of solving the parabolic equation directly,
\begin{align}
    \textbf{X}_{n+1} = & \textbf{X}_n+ \sigma \Delta \textbf{B}_n,  \nonumber\\
    u_{n+1}= &u_n +(\lambda \Psi_{\theta} -\mu^T \nabla \Psi_{\theta})(\textbf{X}_n) \Delta t
+\nabla \Psi_{\theta}(\textbf{X}_n) \Delta \textbf{B}_n.
\end{align}
The original algorithm in \cite{han20} uses a separate DNN to approximate the gradient of $\Psi(\textbf{x})$.
And the loss function for the DNN  $\Psi_{\theta}(\textbf{x})$ approximating the eigen-function and the eigenvalue is then defined by
\begin{equation}
    Loss_{semigroup}(\theta, \lambda) =E_{X_0 \sim \pi_0}[|u_N-\Psi_{\theta}(\textbf{X}_N)|^2].
\end{equation}

\section{Martingale problem formulation of elliptic
PDEs}

In this section, we will present the Martingale problem formulation for the BVPs and eigenvalue problems of elliptic PDEs, and for this purpose,
let us consider a general PDE for the linear elliptic operator $\mathcal{L}$ with $\mu=\mu(\mathbf{x}), \sigma=\sigma(\mathbf{x}) $,
\begin{align}
\mathcal{L}u+V(\mathbf{x}, u, \nabla u)  &  =f(\textbf{x},u),\text{ \ }\mathbf{x}\in D\subset
R^{d},\label{pde}\\
\mbox{or} \quad \mathcal{L}u & =f(\textbf{x}, u)-V(\mathbf{x}, u, \nabla u), \\
\Gamma(u)  &  =g,\text{ \ }\mathbf{x}\in\partial D,\nonumber
\end{align}
where $f(\textbf{x},u)=\lambda u, g=0$ for the case of an eigenvalue problem with an eigenvalue $\lambda$ and eigenfunction $u(x)$, and the boundary operator $\Gamma$ could be Dirichlet, Neumann, or Robin type or
a decay condition will be given at $\infty$ if $D=R^{d}$.
The following shorthand will be used in the rest of the paper
\begin{equation}
   v(\textbf{x})= V(\textbf{x},u(\textbf{x}),\nabla u(\textbf{x})).
   \label{v(x)}
\end{equation}

The vector $\mu = \mu(\mathbf{x}) $, $\sigma_{d\times d} =\sigma_{d\times d}(\mathbf{x})$ can be
associated with the variable drift and diffusion, respectively, of the following stochastic Ito
process $\mathbf{X}_{t}(\omega)$ $\mathbf{\in}$ $R^{d},$ $ \omega\in\Omega$ (random
sample space) with $\mathcal{L}$ as its generator,
\begin{align}
d\mathbf{X}_{t}  &  =\mu (\mathbf{X}_{t}) dt\mathbf{+}\sigma (\mathbf{X}_{t})\mathbf{\cdot}d\mathbf{B}_{t}\label{sde}\\
\mathbf{X}_{t}  &  =\mathbf{x}_{0}\in D, \nonumber
\end{align}
where $\mathbf{B}_{t}=(B_{t}^{1},\cdots,B_{t}^{d})^{\top}\mathbf{\in}R^{d}$ is the
Brownian motion in $R^{d}.$

The transition probability $P(\mathbf{y},t;\mathbf{x},s)$ for the process $\mathbf{X}_{t}$ will satisfy the following Fokker-Planck equation
\begin{equation}
\frac{\partial P(\mathbf{y},t;\mathbf{x},s)}{\partial t}=\mathcal{L}^*_y P(\mathbf{y},t;\mathbf{x},s),
\end{equation}
where the adjoint operator
\begin{equation}
\mathcal{L}^*_y=-\nabla^{\top}_y \mu+ \frac{1}{2}Tr(\nabla_y \nabla_y^{\top}A).
\end{equation}

\bigskip
\noindent ({\bf Robin Problem}) Let us consider the BVP of \eqref{pde} with a Robin type boundary condition
\begin{equation}
\Gamma(u)=\gamma^{\top}\cdot\nabla u+cu=g,\label{RobinBC}%
\end{equation}
where the   vector
\begin{equation}
    \gamma(x) = \frac{1}{2}A\cdot \textbf{n},
\end{equation}
and $\textbf{n}$ is the outward normal at $\textbf{x}\in\partial D$. The Dirichlet BC \ $u=f$
can be considered as a limiting case of $c\rightarrow N,g\rightarrow
Ng,N\rightarrow\infty$ and the Neumann BC as the case $c=0.$

\bigskip With the Martingale problem approach \cite{Varadhan71},
the Martingale problem for the BVP with the third kind boundary condition
(\ref{RobinBC}) can formulated using a reflecting diffusion process
$\mathbf{X}^{ref}$ based on the process $\mathbf{X}$ \eqref{sde} through the
following Skorohod problem.

\bigskip(\textrm{\textbf{Skorohod problem):}} Assume $D$ is a bounded domain
in $R^{d}$ with a $C^{2}$ boundary. Let $\textbf{X}(t)$ be a (continuous) path of \eqref{sde}
in $R^{d}$ with $\textbf{X}(0)\in\bar{D}$. A pair $(\textbf{X}^{ref}(t),L(t))$ is a solution to
the Skorohod problem $S(\textbf{X};D)$ if the following conditions are satisfied:
\begin{enumerate}
\item $\textbf{X}^{ref}$ is a path in $\bar{D}$;

\item ({\bf local time}) $L(t)$ is a non-decreasing function which increases only when $X^{ref}%
\in\partial D$, namely,
\begin{equation}
L(t)=\int_{0}^{t}I_{\partial D}(\textbf{X}^{ref}(s))L(ds),\label{sk1}%
\end{equation}

\item The Skorohod equation holds:
\begin{equation}
S(\textbf{X};D):\qquad\ \textbf{X}^{ref}(t)=X(t)-\int_{0}^{t}\gamma(\textbf{X}^{ref}(s))L(ds).\label{sk2}%
\end{equation}
\end{enumerate}

Here, $L(t)$ is the local time of the reflecting diffusion process, where an
oblique reflection with respect to the direction $\gamma$ at the boundary
$\partial D$ occurs once the process $\textbf{X}(t)$ hits the boundary \cite{schuss15}. The sampling of reflecting process and the computation of local time $L(t)$
can be found in \cite{ding2023}.

As we will only use the reflecting diffusion $\textbf{X}^{ref}(t)$ for the rest of our
discussion, we will keep the same notation
\begin{equation}
    \textbf{X}(t)\leftarrow \textbf{X}^{ref}(t)
\end{equation}
with the
understanding that  it now stands for a reflecting diffusion process within the
closed domain $\bar{D}$. Using the Ito formula \cite{klebaner} for the semi-martingale
$\textbf{X}(t)$ (namely, $\textbf{X}^{ref}(t)$) \cite{hsu84} \cite{Pap1990},%
\[
du(\textbf{X}(t))=
{\displaystyle\sum\limits_{i=1}^{d}}
\frac{\partial u}{\partial x_{i}}(\textbf{X}(t))dX_{i}(t)+\frac{1}{2}
{\displaystyle\sum\limits_{i=1}^{d}}
{\displaystyle\sum\limits_{j=1}^{d}}
a_{ij}(\textbf{X}(t))\frac{\partial^{2}u}{\partial x_{i}\partial x_{j}}(\textbf{X}(t))dt,
\]
with the notation of the generator $\mathcal{L}$, we have for the solution $u(\textbf{x})$ of \eqref{pde} the following differential%
\begin{align*}
& du(\textbf{X}(t)) = \\
& \mathcal{L}u(\textbf{X}(t))dt-\gamma^{\intercal}\cdot\nabla u(u(\textbf{X}_{t}%
))L(dt)+%
{\displaystyle\sum\limits_{i=1}^{d}}
{\displaystyle\sum\limits_{j=1}^{d}}
\sigma_{ij}\frac{\partial u}{\partial x_{i}}(\textbf{X}(t))dB_{i}(t)\\
& =\left[f(\textbf{X}(t),u(\textbf{X}(t))-V(\textbf{X}(t),u(\textbf{X}(t)),\nabla u(\textbf{X}(t)))\right]dt-\left[
g(\textbf{X}(t))-cu(\textbf{X}(t))\right]  L(dt)\\
& +%
{\displaystyle\sum\limits_{i=1}^{d}}
{\displaystyle\sum\limits_{j=1}^{d}}
\sigma_{ij}\frac{\partial u}{\partial x_{i}}(\textbf{X}(t))dB_{i}(t),
\end{align*}
which in turns gives a Martingale $M_{t}^{u}$ defined by
\begin{align}
& M_{t}^{u} \doteq   \\
& u(\textbf{X}_{t})-u(\textbf{X}_{0})-\int_{0}^{t}\left[  f(\textbf{X}_{s}%
,u(\textbf{X}_{s}))-V(\textbf{X}_{s},u(\textbf{X}_{s}),\nabla u(\textbf{X}_{s}))\right]  ds \\
& +\int_{0}^{t}\left[
g(\textbf{X}_{s})-cu(\textbf{X}_{s})\right]  L(ds) =\int_{0}^{t}%
{\displaystyle\sum\limits_{i=1}^{d}}
{\displaystyle\sum\limits_{j=1}^{d}}
\sigma_{ij}\frac{\partial u}{\partial x_{i}}(\textbf{X}_{s})dB_{i}(s), \nonumber
\end{align}
due to the Martingale nature of the Ito integral at the end of the equation above.

\bigskip
\noindent({\bf Dirichlet Problem}) For Dirichlet problem of \eqref{pde} with a boundary condition
\begin{equation}
    \Gamma [u]=u=g, \quad \textbf{x} \in \partial D,
\end{equation}
the underlying diffusion process is
the original diffusion process \eqref{pde}, but killed at the boundary at the first exit time
\begin{equation}
    \tau_{D}=\inf\{t,\textbf{X}_{t}\in\partial D\},
\end{equation}
 and it can be shown that in fact%
\begin{equation}
 \tau_{D}=\inf\{t>0,L(t)>0\},
\end{equation}
and also that $M_{t\wedge\tau_{D}}^{u}$ will be still a Martingale \cite{klebaner}, which
will not involve the integral with respect to local time $L(t)$, i.e.
\begin{equation}
   M_{t\wedge\tau_{D}}^{u}=u(\textbf{X}_{t\wedge\tau_{D}})-u(\textbf{X}_{0})-\int_{0}^{t\wedge
\tau_{D}}\left[  f(\textbf{X}_{s},u(\textbf{X}_{s}))-V(\textbf{X}_{s},u(\textbf{X}_{s}),\nabla u(\textbf{X}_{s}))\right] \label{MartD}
ds.
\end{equation}

For the case of a linear PDE, i.e. $f(\textbf{x}, u)=f(\textbf{x}), V=0$ , by taking expectation of \eqref{MartD} and letting $t \rightarrow \infty$, we will have
\begin{eqnarray}
     0=&E[M^u_{0}]=lim_{t\rightarrow \infty} E[M^u_{t\wedge\tau_{D}}]=E[M^u_{\tau_{D}}] \nonumber\\
     =&E[u(\textbf{X}_{\tau_{D}})-u(\textbf{X}_{0})]-
     \int_{0}^{\tau_{D}}f(\textbf{X}_{s}) ds \nonumber \\
     =& E[g(\textbf{X}_{\tau_{D}})]-u(\textbf{x})-
    E[\int_{0}^{\tau_{D}}f(\textbf{X}_{s}) ds],
\end{eqnarray}
resulting in the well-known Feynman-Kac formuula for the Dirichlet boundary value problem
\begin{equation}
 u(\textbf{x})=  E[g(\textbf{X}_{\tau_{D}})]-E[\int_{0}^{\tau_{D}}f(\textbf{X}_{s}) ds],  \quad \textbf{x} \in D,
\end{equation}
where the diffusion process $\textbf{X}_t, \textbf{X}_0=\textbf{x}$ is defined by \eqref{sde}.

\bigskip The Martingale problem of the BVPs states the equivalence of $M_{t}^{u}$
being a Martingale (i.e. a probabilistic weak form of the BVPs)  and the
classic weak form:
\medskip
For every test function $\phi(\mathbf{x})\in C_{\partial D}^{2}=\{\phi:\phi\in
C^{2}(D)\cap C^{1}(\overline{D}),\left(  \gamma' \cdot\nabla-\mu^{\intercal
}\textbf{n}\right)  \phi=0\}$, $\gamma^{\prime}=\gamma-\alpha$, $\alpha_{j}=%
{\displaystyle\sum\limits_{i=1}^{3}}
\frac{\partial a_{ij} }{\partial x_{i}}$, we have%

\begin{equation}
\begin{aligned}
     \int_{D}u(\textbf{x})\mathcal{L}^{\ast}\phi d\textbf{x}=&\int_{D}\left[
f(\textbf{x},u(\textbf{x}))-V(\textbf{x},u(\textbf{x}),\nabla u(\textbf{x}))\right]  \phi(\textbf{x})d\textbf{x} \\
+&\int_{\partial D}\phi
(\textbf{x})[g(\textbf{x})-cu(\textbf{x})]ds_{x},
\end{aligned}
\end{equation}
where%
\begin{equation}
    \mathcal{L}^{\ast}\phi=\frac{1}{2}Tr(\nabla\nabla^{\intercal}A)\phi
-div(\mu\phi).
\end{equation}
This equivalence has been proven for the Schrodinger operator
$\mathcal{L}u=\frac{1}{2}\Delta u+qu$ for the Neumann problem \cite{hsu84} and the
Robin problem \cite{Pap1990}.

\section{DeepMartNet - a Martingale based neural network}

 In this section, we will propose a DNN method for solving the BVPs and eigenvalue problems for elliptic PDEs using the equivalence between its Martingale problem formulation and classic weak form of the PDEs.

For simplicity of our discussion, let us assume that $s \le t \le \tau_{D}$, by the Martingale property of $M_t\equiv M^u_t$ of \eqref{MartD}, we have
\begin{equation}
E[M_{t}|\mathcal{F}_{s}]=M_{s}, \label{martin}%
\end{equation}
which implies for any measurable set $A\in$ $\mathcal{F}_{s},$

\begin{equation}
E[M_t|A]=M_s=E[M_s|A],
\end{equation}
thus,
\begin{equation}
E[\left(  M_{t}-M_{s}\right) | A ] =0,
\label{Exp1a}%
\end{equation}
i,e,
\begin{equation}
\int_{A}\left(  M_{t}-M_{s}\right)  P(d\omega)=0, \label{Exp}%
\end{equation}
where
\begin{align}
 M_t-M_s
 =&u(\mathbf{X}_t)-u(\mathbf{X}_s)-\int_{s}^{t}\mathcal{L}u(\mathbf{X}%
_{z})dz \nonumber \\
= & u(\mathbf{X}_t)-u(\mathbf{X}_s)-
\int_{s}^{t}(f(z, u(\mathbf{X}_z))-v(\mathbf{X}_z)) dz.
\label{M_ts}
\end{align}

In particular, if we take $A=\Omega\in\mathcal{F}_{s}$ in (\ref{Exp1a}), we have%
\begin{equation}
E[M_{t}-M_{s}]=0, \label{cExp}%
\end{equation}
i.e. the Martingale $M_{t}$ has a constant expectation. However, it should be noted that constant expectation by itself does not mean a Martingale, for this we have the following lemma \cite{cohen15}.
\begin{lemma}
If $E[M_{S}]=E[M_{T}]$ holds for any two stopping time $S\leq T$, then
$M_{t}, t\geq0$ is a Martingale, i.e. $E[M_{t}|\mathcal{F}_{s}]=M_{s}$ for
$s\leq t.$





\end{lemma}

For a
given time interval $[0,T]$, we define a partition
\begin{equation}
    0=t_{0}<t_{1}<\cdots<t_{i}<t_{i+1}<\cdots<t_{N}=T, \label{tmesh}
\end{equation}
and the increment of the $M_t$ over $[t_i, t_{i+k}]$ can be approximated by using a trapezoidal rule for the integral term
\begin{align}
 M_{t_{i+k}}-M_{t_i}
 =&u(\mathbf{X}_{i+k})-u(\mathbf{X}_i)-\int_{t_i}^{t_{i+k}}\mathcal{L}u(\mathbf{X}%
_{z})dz \nonumber \\
\doteq & u(\mathbf{X}_{i+k})-u(\mathbf{X}_i)-
\Delta t\sum_{l=0}^k\omega_l \mathcal{L}u(\mathbf{X}_{i+l}) \nonumber \\
=& u(\mathbf{X}_{i+k})-u(\mathbf{X}_i)-
\Delta t \sum_{l=0}^k\omega_l (f(\mathbf{X}_{i+l}, u(\mathbf{X}_{i+l}))-v(\mathbf{X}_{i+l})),
 \label{DeltaM}
\end{align}
where for $k\ge 1$ $\omega_0=\omega_k=\frac{1}{2},\quad \omega_l=1, \quad 2\le l \le k-1$ and for $k=0$, $\omega_0=1$.

Adding back the exit time $\tau_D$, we note that
\[
M_{t_{i+k}\wedge\tau_D}-M_{t_i\wedge\tau_D}=
u(\mathbf{X}_{t_{i+k}\wedge\tau_D})-u(\mathbf{X}_{t_i\wedge\tau_D})-\int_{t_i\wedge\tau_D}^{t_{i+k}\wedge\tau_D}\mathcal{L}u(\mathbf{X}%
_{z})dz=0
\]
if both $t_{i+k}, t_i\ge \tau_D$.



\begin{remark}
We could define a different generator $\mathcal{L}$ by not including ${\mu}%
^{\top}\nabla$ in (\ref{gen}), then the Martingale in (\ref{Mart1}) will be
changed to
\begin{equation}
M_{t}^{\ast}=u(\mathbf{X}_{t})-u(\mathbf{x}_{0})-\int_{0}^{t}(\lambda u(\mathbf{X}_{s})
-\mu^{\top}(\mathbf{X}_{s})\nabla u(\mathbf{X}_{s})-v(\mathbf{X}_{s}))ds,
\label{Mart1}%
\end{equation}
where the process $\mathbf{X}_{t}$ is given simply by $d\mathbf{X}_{t}=\sigma
\mathbf{\cdot}d\mathbf{B}_{t},$ instead.
\end{remark}%

\begin{itemize}
\item \bigskip {\bf DeepMartNet for Dirichlet BVPs}
\end{itemize}

Let $u_{\theta}(\mathbf{x})$ be a neural network approximating the
BVP solution with $\theta$ denoting all the weight and bias parameters of a DNN. For a
given time interval $[0,T]$, we define a partition%
\[
0=t_{0}<t_{1}<\cdots<t_{i}<t_{i+1}<\cdots<t_{N}=T,
\]
and $M$-discrete realizations
\begin{equation}
\Omega^{\prime}=\{\omega_{m}\}_{m=1}^{M}\subset\Omega\label{M-sample}%
\end{equation}
of the Ito process using the Euler-Maruyama scheme with $M$-realizations of the
Brownian motions $\mathbf{B}_{i}^{(m)}$, $0\leq m\leq M$, $\mathbf{\ }$%

\[
\mathbf{X}_{i}^{(m)}(\omega_{m})\sim X(t_{i},\omega_{m}),0\leq i\leq N,
\]
where%
\begin{align}
\mathbf{X}_{i+1}^{(m)}  &  =\mathbf{X}_{i}^{(m)}+\mu\mathbf{(X}_{i}%
^{(m)}\mathbf{)}\Delta t_{i}\mathbf{+}\sigma\mathbf{\mathbf{(X}_{i}%
^{(m)}\mathbf{)}\cdot}\Delta\mathbf{B}_{i}^{(m)},\label{euler}\\
\mathbf{X}_{0}^{(m)}  &  =\mathbf{x}_{0}
\end{align}
where $\Delta t_{i}=t_{i+1}-t_{i},$%
\[
\Delta\mathbf{B}_{i}^{(m)}=\mathbf{B}_{i+1}^{(m)}-\mathbf{B}_{i}^{(m)}.
\]

\bigskip

We will build the loss function $Loss(\theta)$ for
neural network $u_{\theta}(\mathbf{x})$ approximation of the BVP solution using the
Martingale property (\ref{Exp}) and the M-realization of the Ito diffusion
(\ref{sde}).

For each $t_{i}$, we randomly take a subset of $A_{i}\subset\Omega^{\prime}$
with a uniform sampling (without replacement), corresponding to the mini-batch
in computing the stochastic gradient for the empirical training loss. Assuming that the mini-batch in each $A_{i}$ is large enough such that
$\{\mathbf{X}_{i+1}^{(m)}\}$ and $\{\mathbf{X}_{i}^{(m)}\}$, $\omega_m \in A_{i} $ sample the distribution $P(t_{i+1},.,t_i,A), P(t_i,.,0,x_0)$ well, respectively, then equation \eqref{ExpA} gives
the following approximate identity for the solution $u_{\theta}(\mathbf{X}_{t})$
using the $A_{i}-$ensemble average,%

As $E[M_{t_{i+k}}^{u_{\theta}}-M_{t_i}^{u_{\theta}} ] \approx 0$, for a randomly selected $A_i \in \Omega^{\prime} = \mathcal{F}_{t_i}$ ( mini-batches), we can set

 \begin{align}
    Loss_{mart}(\theta)  & =\frac{1}{N}
   \sum_{i=0}^{N-1} (M_{t_{i+k}\wedge\tau_D}^{u_{\theta}}-M_{t_i\wedge\tau_D}^{u_{\theta}})^2 \nonumber \\
    & = \frac{1}{N}\sum_{i=0}^{N-1}\frac{I({t_i \le \tau_D})}{|A_i|^2}\left(\sum_{m=1}^{|A_i|}
      u_{\theta}(\mathbf{X}_{i+k}^{(m)})-u_{\theta}(\mathbf{X}_{i}^{(m)})- \right.  \nonumber\\
& \left.  \Delta t\sum_{l=0}^k\omega_l (f(\mathbf{X}_{i+l}^{(m)}, u_{\theta}(\mathbf{X}_{i+l}^{(m)}))
-v_{\theta}(\mathbf{X}_{i+l}^{(m)}))  \right)^{2},
\label{loss_Mart_bvp}
\end{align}
where $v_{\theta}(\textbf{x})$ is defined similarly as in \eqref{v(x)}. Refer to Remark \ref{remark-batch} for the discussion of the size of the mini-batch $|A_i|$.

Now, we define the total loss for the boundary value problem as
\begin{equation}
   Loss_{total-bvp}(\theta)  =Loss_{mart}(\theta)+ \alpha_{bdry} Loss_{bdry}(\theta),
   \label{totallossBvp}
\end{equation}
and $\alpha_{bdry}$ is a hyper-parameter and $Loss_{bdry}(\theta)=||u_{\theta}(\textbf{x})-g||_2^2$, which can be approximated by evaluations at the boundary.

 \medskip
\textbf{DeepMartNet solution}- $u_{\theta^*}(x)$, here
\begin{equation}
    \theta^*=argmin Loss_{total-bvp}(\theta).
\end{equation}

\begin{itemize}
\item \bigskip {\bf DeepMartNet for Dirichlet eigenvalue problems}
\end{itemize}

For the eigenvalue problem
\begin{align}
\mathcal{L}u+V(\mathbf{x}, u, \nabla u)  &  =\lambda u,\text{ \ }\mathbf{x}\in D\subset
R^{d},\label{pde_eig}\\
\Gamma(u)=u  &  =0,\text{ \ }\mathbf{x}\in\partial D,\nonumber
\end{align}
and the Martingale loss becomes
\begin{align}
    Loss_{mart}(\lambda, \theta)   =\frac{1}{N}\sum
_{i=0}^{N-1}  & \frac{1}{|A_{i}|^2}\left(\sum_{m=1}^{|A_{i}|}  u_{\theta
}(\mathbf{X}_{i+k}^{(m)})-u_{\theta}(\mathbf{X}_{i}^{(m)})- \right. \nonumber\\
& \left. \Delta t\sum_{l=0}^k\omega_l (\lambda u_{\theta}(\mathbf{X}_{i+l}^{(m)})
-v_{\theta}(\mathbf{X}_{i+l}^{(m)}))
\right) ^{2}
\label{loss_Mart_eig}
\end{align}
and in the case of a bounded domain, the boundary loss $Loss_{bdry}(\theta)$ will be added for the homogeneous boundary condition $g=0$, i.e., $Loss_{bdry}(\theta)=||u_{\theta}(\textbf{x})||_2^2$. For the decay condition at infinite, a molifier will be used to enforce explicitly the decay condition there (see \eqref{mollifier} ).

Also, in order to prevent the DNN eigenfunction going to a zero solution, we introduce a simple normalization term using $l_p$ (p=1, 2) norm of the solution at some randomly selected location
\begin{equation}
    {Loss}_{normal}(\theta) =  \left( \frac{1}{m}\sum_{i = 1}^m |u_{\theta}({\textbf{x}_i})|^p - c \right)^2,
    \label{lossNormal}
\end{equation}
where  $\textbf{x}_i$ are m arbitrarily selected fixed points and $c$ is a nonzero constant.

Finally, we have the total loss for the eigenvalue problem as
\begin{equation}
    Loss_{total-eig}(\lambda, \theta)= Loss_{mart}(\lambda, \theta)+ \alpha_{bdry}Loss_{bdry}(\theta)+ \alpha_{normal} {Loss}_{normal}(\theta),
    \label{totallossEig}
\end{equation}
where $\alpha_{bdry}$ and $\alpha_{normal}$ are hyper-parameters.

\medskip
\textbf{DeepMartNet eigen-problem solution }- $(\lambda^*, u_{\theta^*}(x))$,  here
\begin{equation}
    (\lambda^*, \theta^*)=argmin Loss_{total-eig}(\lambda, \theta).
\end{equation}
\begin{remark} ({\bf Mini-batch in SGD training and Martingale property})
Due to the equivalence between (\ref{Exp}) and (\ref{martin}), the loss
function defined above ensures that $M_{t}^{u_{\theta}}$ of (\ref{MartD}) for $u_{\theta
}(\mathbf{x})$ will be a Martingale approximately if the mini-batch $A_{i}$
explores all subsets of the sample space $\Omega^{\prime}$ during the SGD optimization process of
the training, and the sample size
$M=|\Omega^{\prime}|$ $\rightarrow$ $\infty,$ the time step $\max|\Delta
t_{i}|\rightarrow0,$ and the training converges (see Fig. \ref{fig:mart-sgd}).

Also, if we take $A_{i}=\Omega^{\prime}$ for all $i,$ there will be no stochasticity
in the gradient calculation for the loss function, we will have a traditional
full gradient descent method and the full Martingale property for  $u_{\theta
}(\mathbf{x})$  is not enforced either. Therefore, the mini-batch procedure in DNN SGD optimization
corresponds perfectly with the Martingale definition (\ref{martin}).

In summary, the Martingale property implies that for any measurable set $A\in$ $\mathcal{F}_{s},$ we require that
    \begin{equation}
E[M_{t}|A]=M_{s},
\end{equation}
which then provides a native mechanism for the mini-batch in the SGD.  Therefore, the Martingale based DNN is an ideal fit for deep learning of high-dimensional PDEs.
\begin{figure}[htb]
    \centering
    \includegraphics[width=0.8\textwidth]{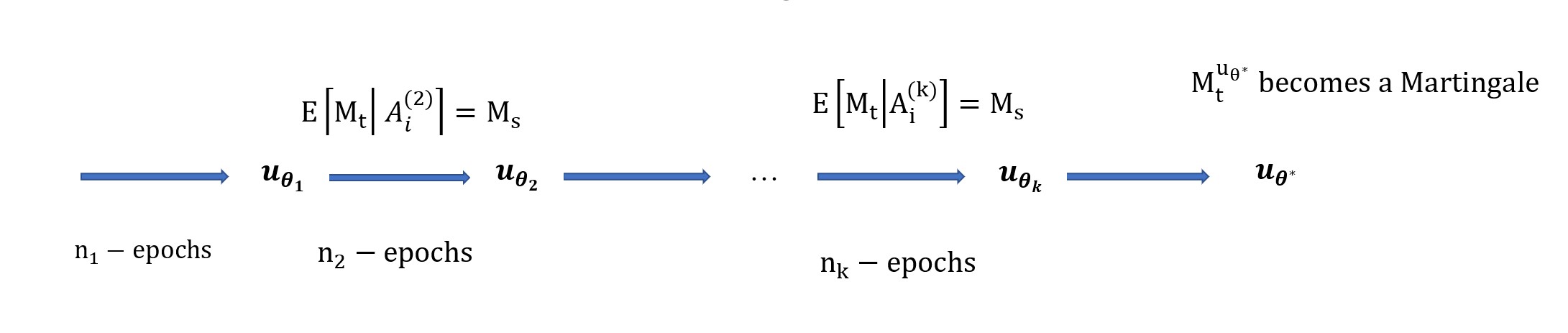}
    \caption{DeepMartNet training and Martingale property}
    \label{fig:mart-sgd}
\end{figure}

\end{remark}

\begin{remark}
\label{remark-batch}
    ({\bf Size of mini-batch $A_i$})
The loss function of the DeepMartNet is based on the fact that $\int_{A_i}\left(  M^u_{t}-M^u_{s}\right)  P(d\omega)=0 $ for the exact solution $u(\textbf{x})$, where the expectation will be computed by an ensemble average of selected paths from the $M$-paths. In theory, using the transitional probability, using left end point quadrature for the integral over [s,t] in \eqref{M_ts} with $|t-s|<<1$, \eqref{Exp} can be rewritten as
\begin{align}
& 0=  E[\left(  M^u_{t} -M^u_{s}\right)  | A_i ]  \nonumber \\
 \doteq & \int_{\textbf{z} \in B }\int_{\textbf{y}\in R^d} \big[ u(\textbf{y})-u(\textbf{z})-(f(\textbf{z}, u(\textbf{z})) -v(\textbf{z}))(t-s) \big] P(t,\textbf{y},s,\textbf{z})p(s,\textbf{z},0,\textbf{x}_0)dzdy \nonumber \\
  = &\int_{\textbf{y}\in R^d} u(\textbf{y}) P(t,\textbf{y},s,A_i)dy  -\int_{\textbf{z} \in B } \big[u(\textbf{z})+(f(\textbf{z}, u(\textbf{z})) -v(\textbf{z}))(t-s)\big] P(s,\textbf{z},0,\textbf{x}_0) dz \nonumber \\
 \sim & \frac{1}{|A_i|} \sum_{m=1}^{|A_i|} \left( u(\textbf{X}_t^{(m)})-
 \big[u(\textbf{X}_s^{(m)})+(f(\textbf{X}_s^{(m)}, u(\textbf{X}_s^{(m)})) -v(\textbf{X}_s^{(m)}))(t-s)\big]
 \right),  
\label{ExpA}
\end{align}
where $B=\mathbf{X}^{-1}_{s}(A_i)$ and $\textbf{X}_t^{(m)},  1 \le m \le M$ are the M-sample paths of the diffusion process $\textbf{X}_t$.

 Therefore, the size of mini-batch $|A_i|$ should be large enough to give an accurate sampling of the  continuous distribution
 $P(t,\textbf{y},s,A), y \in R^d$ and $P(s,\textbf{z},0,\textbf{x}_0), z \in B$, so in our simulation we select a sufficient large $M$ and set the size of mini-batch $A_i$ in the following range,
 \begin{equation}
     M/m_1 \le |A_i| \le M/m_2,
     \label{bsize}
 \end{equation}
 where $M$ is the total number of paths used and $m_1>m_2$ are hyper-parameters of the training.

 The set $A_i$ can be the same for each $0\le i \le N-1$ per epoch or can be randomly selected as the subset of $\Omega'$ depending how much stochasticity is put into the calculation of the stochastic gradient in the SGD optimizations.

For a low memory implementation of the DeepMartNet, at any time we can just generate enough $|A_i|$ number of paths to be used for some epochs of training, and then regenerate them again without first generating a large number of paths upfront.

\end{remark}

 \section{Numerical Results}

 The numerical parameters for the DeepMartNet consist of
 \begin{itemize}
     \item M - number of total paths
     \item T - terminal time of the paths
     \item N - number of time steps over [0, t]
     \item $\Delta t = T/N$
     \item $M_b=|A_i|$ size of mini-batches of paths selected according to \eqref{bsize}.
     \item Size of networks
 \end{itemize}
 The training is carried out on a Nvidia Superpod one GPU node A100. The Optmizer is Adamax \cite{adam}.

 \subsection{Dirichet BVPs of the Poisson-Boltzmann equation }
We will first apply the DeepMartNet to solve the Dirichlet BVP of the Poisson-Boltzmann equation (PBE) arising from solvation of biomolecules in inoic solvents \cite{caibook}.
 \begin{equation}
\begin{cases}
    \Delta u(\textbf{x}) + c u(\textbf{x}) = f(\textbf{x}),& \textbf{x} \in D \\
    u(\textbf{x}) = g(\textbf{x}),& \textbf{x} \in \partial D
\end{cases}
\end{equation}
where $c<0$  ($c=-1$ in the numerical tests) with an high-dimensional solution given by
\begin{equation} \label{PBEsol}
    u(\textbf{x}) = \sum_{i=1}^d \cos (\omega x_i), \quad \omega=2.
\end{equation}

In this case, the generator for the stochastic process is $\mathcal{L}=\frac{1}{2}\Delta$, so the corresponding diffusion is simply the Brownian motion $\textbf{B}(t)$.
For the $M-$ Brownian paths $\textbf{B}^{(j)}, j=1, \cdots, M$ originating from $x_0$, the Martingale loss \eqref{loss_Mart_bvp} becomes
\begin{align}
Loss_{\text{mart}}( \theta ) := &\frac{1}{N}\sum_{i=0}^{N-1}\frac{1}{|A_i|^2} \left( \sum _{j=1}^{|A_i|} u_{\theta }( \textbf{B}_{t_{i+1}}^{( j)}) -u_{\theta }( \textbf{B}_{t_{i}}^{( j)}) \right.\nonumber \\
-& \left. \frac{1}{2} \left( f ( \textbf{B}_{t_{i}}^{( j)} )  - cu(\textbf{B}_{t_i}^{(j)}) \right)
\mathbb{I} ( t_{i} \leq \tau _{D }^{( j)} ) \Delta t \Big) \right)^{2}.
\end{align}
Meanwhile, we could use the Feynman-Kac formula \cite{klebaner} to compute the solution at the point $\textbf{x}_0$ by
\begin{equation}
  u(\textbf{x}_0)  \approx \frac{1}{M}\sum _{j=1}^{M}\left( g\left( \textbf{B}_{\tau _{D}}^{( j)}\right) \mathrm{e}^{\frac{c\tau_{D}}{2}} +\frac{1}{2}\sum_{i=0}^{N-1} f\left( \textbf{B}_{t_{i}}^{( j)}\right)\mathrm{e}^{\frac{c t_i}{2}}\mathbb{I}\left( t_{i} \leq \tau _{D }^{( j)}\right)   \Delta t\right),
\end{equation}
and also define a point solution loss (or an integral identity for the solution for PDE with quasi-linear term $V$ and $f$ in \eqref{pde} ), termed Feynman-Kac loss, which is added to the total loss \eqref{totallossBvp}
\begin{equation}
    Loss_{\text{F-K}}( \theta )  = (u_{\theta}(\textbf{x}_0) - u(\textbf{x}_0))^2,
\end{equation}
and a boundary loss is approximated as
\begin{equation}
    Loss_{bdry}( \theta )  =||u_{\theta}-g||^2_2 \approx \frac{1}{N_{bdry}}\sum^{N_{bdry}}_{k=1} |u_{\theta}(\textbf{x}_k)-g(\textbf{x}_k)|^2,
\end{equation}
where uniformly sampled boundary points $\textbf{x}_k, 1 \le k \le N_{bd} $ are used to compute the boundary integral.
Therefore, the total loss for the boundary value problem of the PB equation is
\begin{equation}
    Loss_{bvp}( \theta ) =Loss_{\text{mart}}( \theta ) + \alpha_{\text{F-K}} Loss_{\text{F-K}}( \theta )  + \alpha_{bdry} Loss_{bdry}( \theta ),
    \label{loss_bvp}
\end{equation}
where the penalty parameter $\alpha_{\text{F-K}}$ ranges from 10 to 1000 and $\alpha_{bdry}$ ranges from 1000 to 10,000.
An Adamax optimizer with learning rate 0.05 is applied for training  and $\alpha _{\text{F-K}} = 10,  \alpha _{\text{bdry}} = 10^3$ are taken for the following numerical tests.
\begin{itemize}
    \item  {\bf Test 1: PBE in a d=20 dimensional cube $[-1,1]^d$.} In this test, we solve the PBE in a 20 dimensional cube with the Dirichlet boundary condition given by the exact solution \eqref{PBEsol}. The total number of paths is $M = 100,000$ starting from origin, and mini-batch
size $M_b =|A_i|= 1000; \Delta t = 0.01$ and T=9.  A fully connected network with layers (20, 64, 32, 1) with first hidden layer Tanh and second hidden layer GeLU activation function is used.  A mini-batch of $N_\text{bdry} = 2000$ points $\textbf{x}_k$ are uniformly sampled points on the boundary for every epoch of training.

Fig. \ref{fig:PB_20d} show the learned solution along diagonal of the cube (top left) as well as along the first coordinate (top right) and the history of the loss (bottom left) and relative error L2 (computed by Monte Carlo sample) (bottom right). The training takes about less than 5 minutes.
\begin{figure}[htb]
 \centering
 \includegraphics[width=0.35\textwidth]{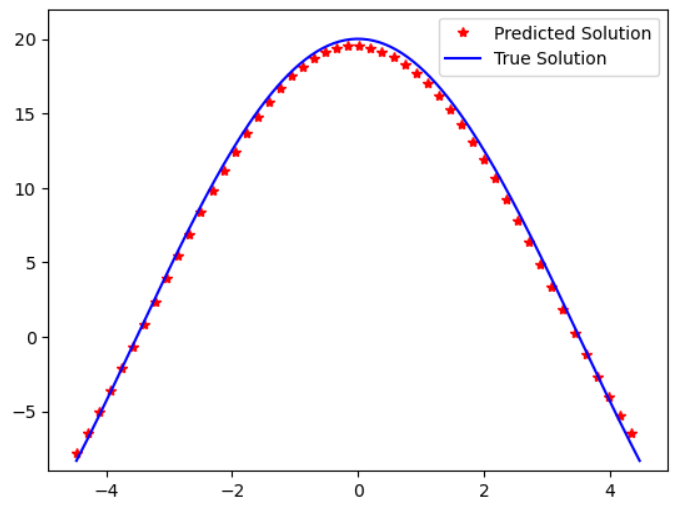}
        \includegraphics[width=0.35\textwidth]{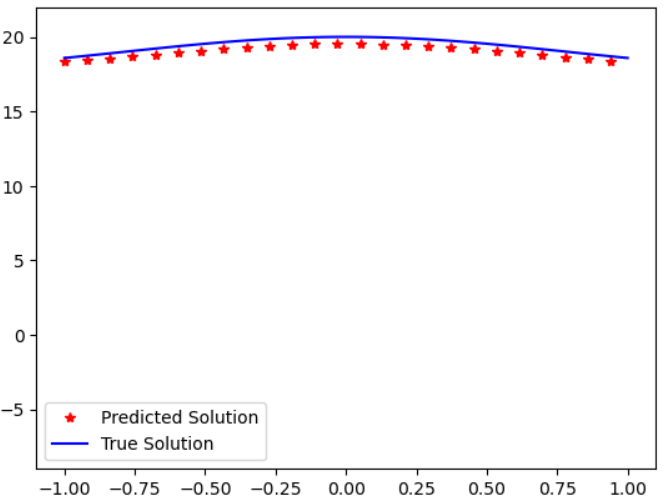}
    \includegraphics[width=0.35\textwidth]{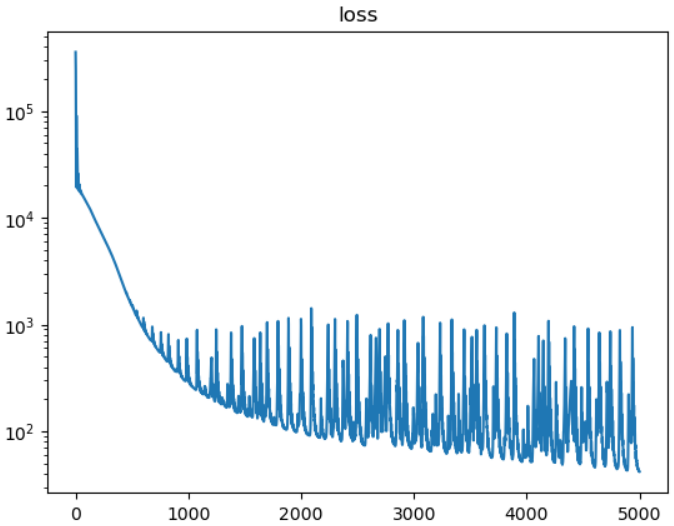}
    \includegraphics[width=0.35\textwidth]{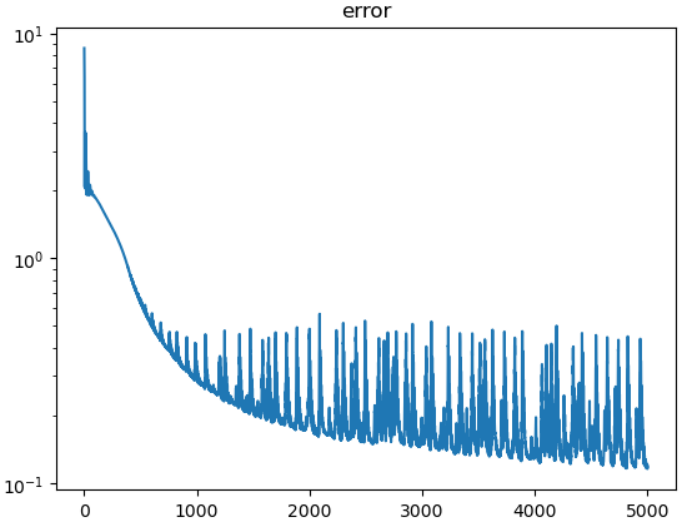}
    \caption{ DeepMartNet solution of PBE in $D = [-1,1]^{d}, d=20$.
    (Upper left): true and predicted value of $u$ along the diagonal of the unite cube; (Upper right): true and predicted value of $u$ along the first coordinate axis of $\mathbb{R}^d$. (Lower left): The loss $L$ history; (lower right): The history of relative error $L^2$ over the cube.}
    \label{fig:PB_20d}
\end{figure}

\item{\bf Test 2. Effect of path starting point $\textbf{x}_0$}. In previous test, the DeepMartNet uses all diffusion paths originating from a fixed point $\textbf{x}_0$ to explore the solution domain. In this test, we consider paths starting from different $\textbf{x}_0=(l,0,\cdots, 0), l=0.1,0.3,0.7$ to investigate the effect of different starting point $\textbf{x}_0$ on the accuracy of the DeepMartNet.
A fully connected network with layers (20, 64, 32, 1) with first hidden layer Tanh and second hidden layer GeLU activation function will be used as in Test 1.  The total number of paths $M = 100,000$, and for each epoch, we randomly choose $M_{b} = 1000$ from the paths; the time step of the paths is $\Delta t = 0.01$.

Fig. \ref{fig:PB_20d3x0} shows that the three different choices of $\textbf{x}_0$ produce similar numerical results with the same numerical parameters as in Test 1.
    \begin{figure}[htbp]
    \centering
       \includegraphics[width=0.35\textwidth]{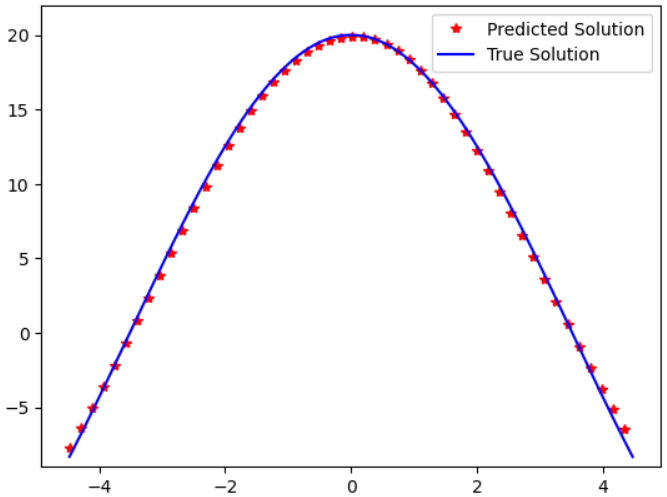}
    \includegraphics[width=0.35\textwidth]{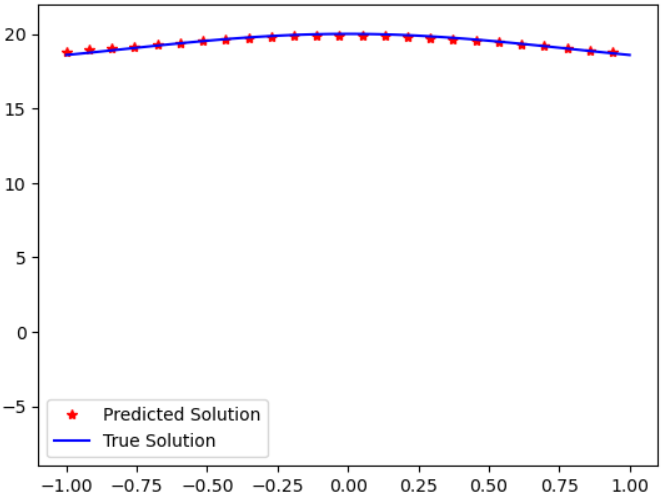}
  \includegraphics[width=0.35\textwidth]{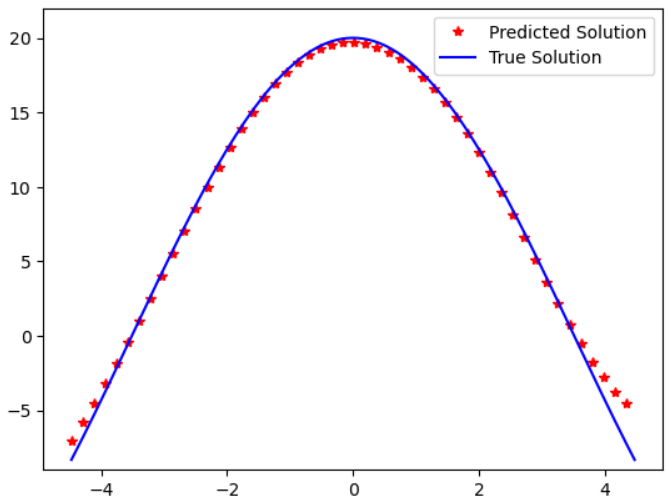}
    \includegraphics[width=0.35\textwidth]{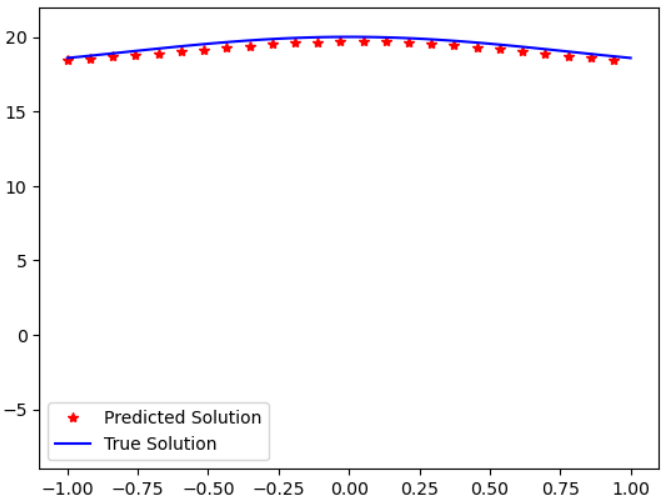}
    \includegraphics[width=0.35\textwidth]{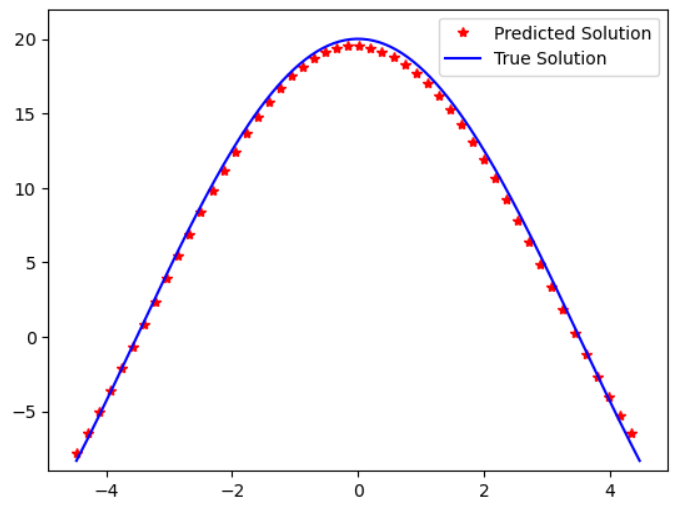}
    \includegraphics[width=0.35\textwidth]{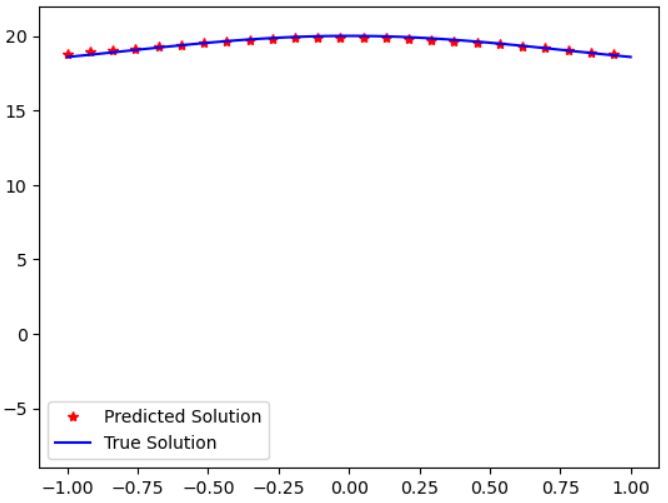}
    \caption{
    DeepMartNet for PBE in $D = [-1,1]^{d}, d=20$ with 3 different starting points for the paths.
    From up to down: $\textbf{x}_0=(l,0,\cdots, 0), l=0.1,0.3,0.7.$  (left):$u(x\boldsymbol{e}_1)$ where $\boldsymbol{e}_1 = d^{-1/2}(1,1\cdots,1)$,(Right): $u(x\boldsymbol{e})$ where $\boldsymbol{e} = (1,0\cdots,0)$.}
    \label{fig:PB_20d3x0}
\end{figure}

Moreover, the DeepMartNet can use diffusion paths starting from different initial position $\textbf{x}_0$ in training the DNN as long as the mini-batch of the paths $A_i$ for a given epoch corresponds to paths originating from a common initial point $\textbf{x}_0$. In Fig. \ref{3start_point}, we compare the numerical result using 120,000 total paths starting from $\textbf{x}_0=(0.3,0,\cdots, 0,0)$ with that with three choices of $\textbf{x}_0$ as in Fig. \ref{fig:PB_20d3x0} with 40,000 paths for each $\textbf{x}_0$. In both cases, for each epoch, we randomly choose mini-batch size $M_b = 1000$ from the paths; the time step of the paths is $\Delta t = 0.01$. An Adamax optimizer with learning rate 0.05 is applied for training. $\alpha _{\text{F-K}} = 10, \alpha_{\text{bdry}} = 10^3$.
The DeepMartNet produce similar results for these two cases with same other numerical parameters as in Test 1.

\begin{figure}
    \centering
    \includegraphics[width=0.45\textwidth]{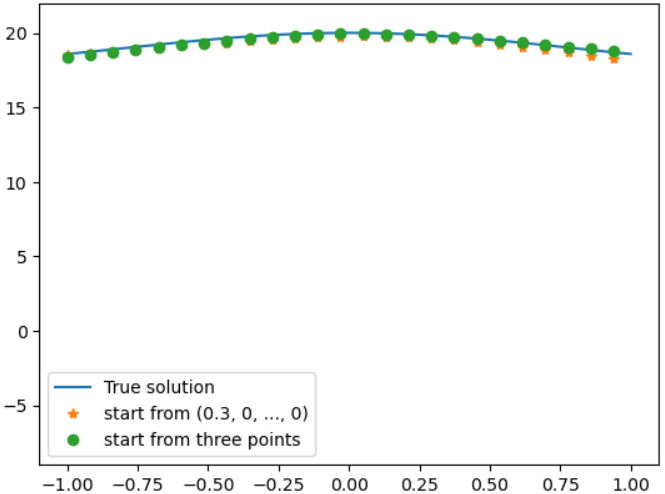}
    \includegraphics[width=0.45\textwidth]{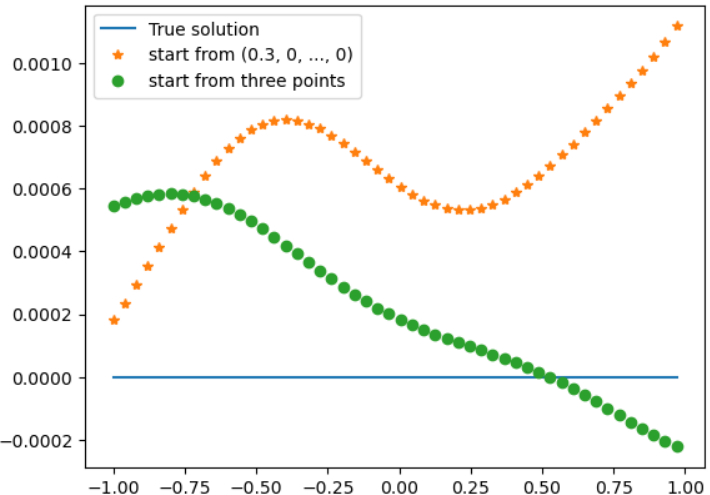}
    \includegraphics[width=0.45\textwidth]{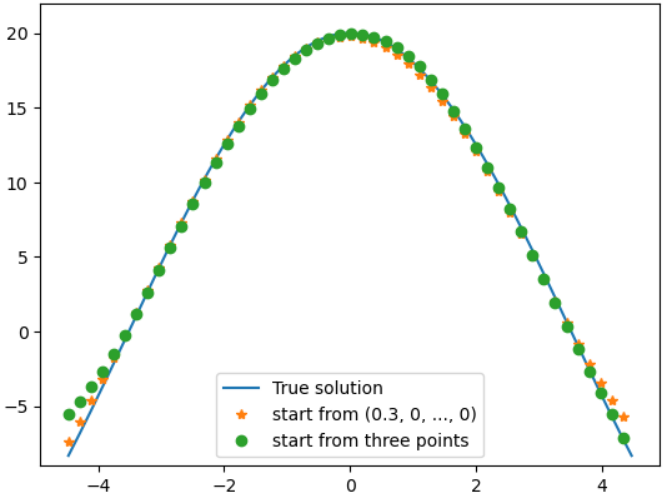}
    \includegraphics[width=0.45\textwidth]{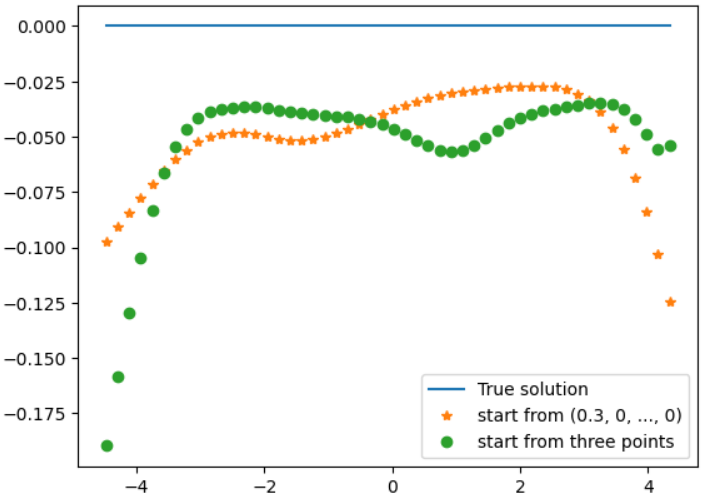}
    \caption{Comparison of DeepMartNet with paths from one starting point and from 3 starting points with same number of total paths in both cases.
     Upper left: $u(x\boldsymbol{e})$ where $\boldsymbol{e} = (1,0\cdots,0)$; Upper right: the relative error $|u-u_\text{true}|/\|u_\text{true}\|_{\infty}$ for the upper left plot; Lower left: $u(x\boldsymbol{e}_1)$ where $\boldsymbol{e}_1 = d^{-1/2}(1,1\cdots,1)$; Lower right: the error $\|u_\text{true}\|_{\infty}$ for the lower left plot.}
    \label{3start_point}
\end{figure}

\item {\bf Test 3: PBE in a d=100 dimensional unit ball.} In this test, we solve the PBE in a 100 dimensional space, We set $T=0.25$ and the time step of the paths is $\Delta t = 0.005$, and
the total number of paths is $M = 100,000$, and for each epoch, we randomly choose $M_b = 1000$ from the paths. A fully connected network with layers (100, 128, 32, 1) with first hidden layer Tanh and second hidden layer GeLU activation function is used.

\begin{figure}[htbp]
    \centering
    \includegraphics[width=0.35\textwidth]{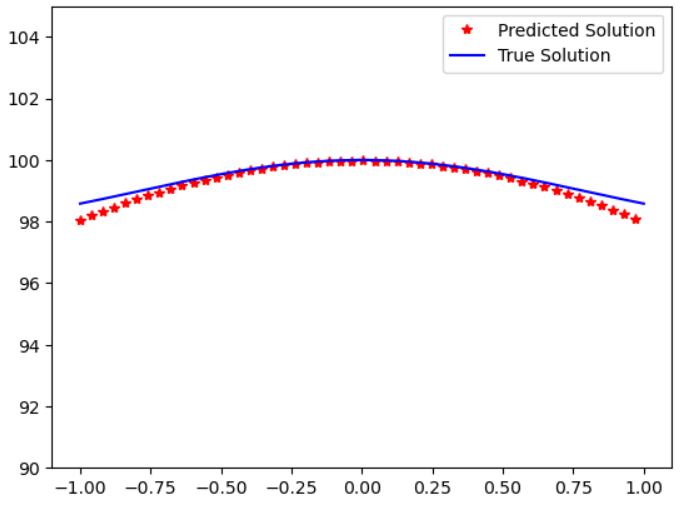}
    \includegraphics[width=0.35\textwidth]{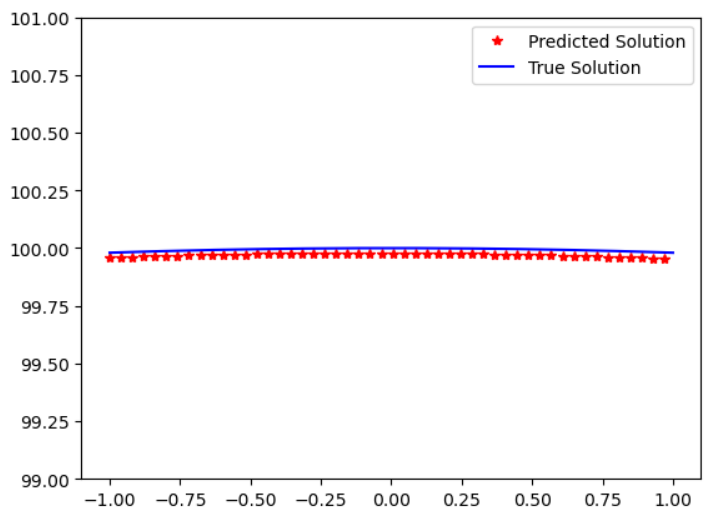}

    \includegraphics[width=0.35\textwidth]{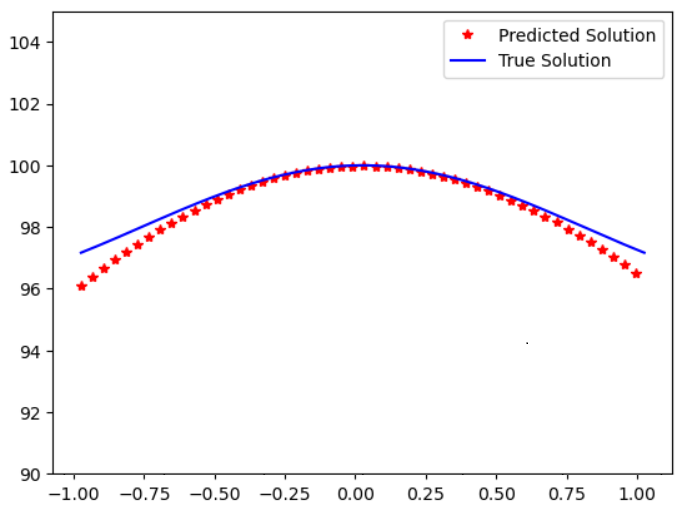}
    \includegraphics[width=0.35\textwidth]{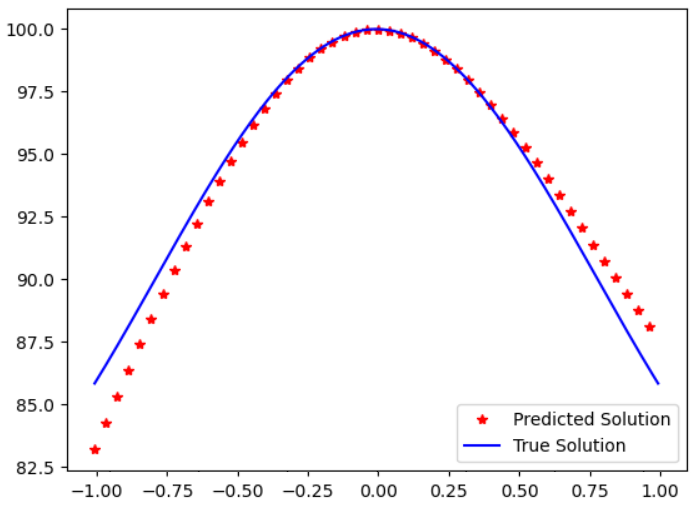}

    \includegraphics[width=0.35\textwidth]{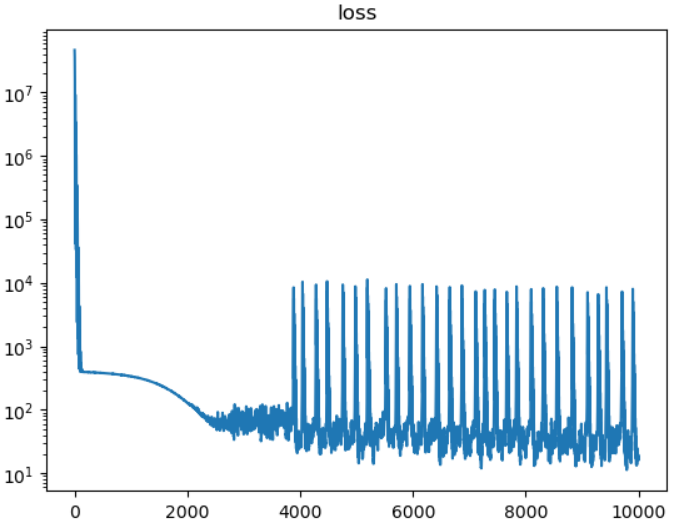}
    \includegraphics[width=0.35\textwidth]{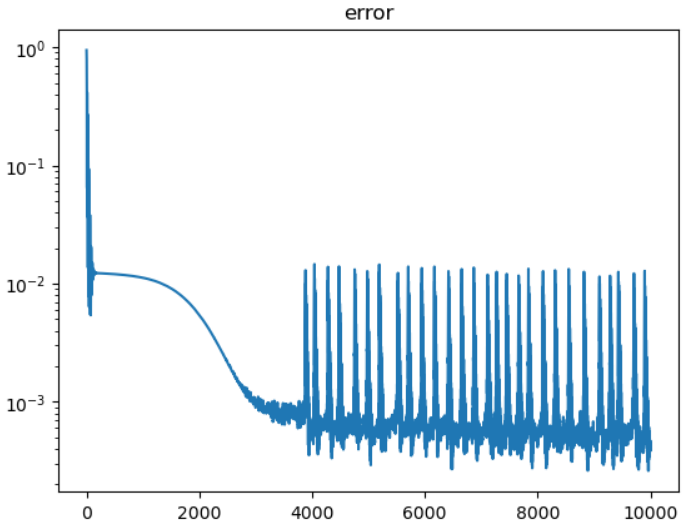}
    \caption{
    DeepMartNet solution of the PBE in $D = \{\textbf{x}\in \mathbb{R}^d, |\textbf{x}| < 1\}$, d=100.
    (Upper left): true and predicted value of $u$ at the diagonal of $\mathbb{R}^d$, i.e. $x$ versus $u(x\boldsymbol{e})$ where $\boldsymbol{e} = d^{-1/2}(1,1\cdots,1)$; (Upper right): true and predicted value of $u$ at the first coordinate axis of $\mathbb{R}^d$, i.e. $x$ versus $u(x\boldsymbol{e_1})$ where $\boldsymbol{e} = d^{-1/2}(1,0\cdots,0)$.
   ( Middle left): $u(x\boldsymbol{e}')$ where $\boldsymbol{e}' = 2^{-1/2}(1,1,0,0,\cdots,0)$
    (Middle right): $u(x\boldsymbol{e}'')$ where $\boldsymbol{e}'' = 10^{-1/2}(\underset{ten 1's}{1,\cdots,1}, 0, 0,\cdots, 0)$
    (Lower left): The loss $L$ versus the number of epoch; (lower right): The relative error $L^2$ error versus the number of epoch.}
    \label{fig:PB_100d}
\end{figure}

 Fig. \ref{fig:PB_100d} show the learned solution (Upper left): true and predicted value of $u$ at the diagonal of $\mathbb{R}^d$, i.e. $x$ versus $u(x\boldsymbol{e})$ where $\boldsymbol{e} = d^{-1/2}(1,1\cdots,1)$; (Upper right): true and predicted value of $u$ at the first coordinate axis of $\mathbb{R}^d$, i.e. $x$ versus $u(x\boldsymbol{e_1})$ where $\boldsymbol{e} = d^{-1/2}(1,0\cdots,0)$.
   ( Middle left): $u(x\boldsymbol{e}')$ where $\boldsymbol{e}' = 2^{-1/2}(1,1,0,0,\cdots,0)$
    (Middle right): $u(x\boldsymbol{e}'')$ where $\boldsymbol{e}'' = 10^{-1/2}(\underset{ten 1's}{1,\cdots,1}, 0, 0,\cdots, 0)$
    (Lower left): The loss $L$ versus the number of epoch; (lower right): The relative error $L^2$ error versus the number of epoch.  The training takes about less than 30 minutes.


\item{\bf Test 4: nonlinear PBE in a d=10 dimensional unit ball.}
For a 1:1 symmetric two species ionic solvent, the electrostatic potential based on the Debye-Huckel theory \cite{caibook} is given by a nonlinear PBE before linearization, and we will consider the following model nonlinear PBE to test the capability of the DeepMartNet for solving nonlinear PDEs,
\begin{equation} \label{nonlinPBsinh}
    \begin{cases}
    -\triangle u +\sinh u = f & x \in D, \\
    u = g & x \in \partial D,
    \end{cases}
\end{equation}
where
\begin{equation}
D = \{x \in {\mathbb R}^{d}: \|\textbf{x}\|_{2} \leq L\}, L=1.
\end{equation}
We consider the case of a true solution as
\begin{equation}
u = \alpha \sum_{i=1}^d x_i ^2,  \quad \alpha=2,
\end{equation}
which gives the boundary data and right hand side of the PBE as
$$
g \equiv \alpha L^2,
$$
and
$$
f = -2\alpha d + \sinh \left( \alpha \sum_{i=1}^d x_i^2 \right).
$$

This time, we use the loss as
\begin{equation}
Loss( \theta ) = Loss_{\text{mart}}( \theta ) + \alpha_{\text{bdry}}Loss_{\text{bdry}}( \theta ),
\end{equation}
where $\alpha_{\text{bdry}} = 10^{-4}$,  and the Martingale loss  \eqref{loss_Mart_bvp} is
\begin{align}
   & Loss_{\text{mart}}( \theta ) = \frac{1}{N}\sum_{i=0}^{N-1}\frac{\mathbb{I}\left( t_{i} \leq \tau _{\partial D }^{( j)}\right)}{M_b^2}  \\
 &\cdot \left(\sum _{j=1}^{M_b} u_{\theta }\left( W_{t_{i+1}}^{( j)}\right) -u_{\theta }\left( W_{t_{i}}^{( j)}\right)   -\frac{\Delta t}{2}\left(\sinh u(W_{t_i}^{(j)})-f\left( W_{t_{i}}^{( j)}\right) \right)  \right)^{2} , \nonumber
\end{align}

and the boundary condition loss is again
\begin{equation} \label{lbdry''}
Loss_{\text{bdry}}( \theta ) =\frac{1}{N_{\text{bdry}}}\sum _{k=0}^{N_{\text{bdry}}}( u_{\theta }( \textbf{x}_{k}) -g(\textbf{x}_{k}))^{2},
\end{equation}
and a mini-batch of $N_\text{bdry} = 2000$ points $\textbf{x}_k$ are uniformly sampled points on the boundary for every epoch of training. And in this case, we do not have a Feynman-Kac loss at the starting point $\textbf{x}_0$.

The numerical results in Fig. \ref{fig:NonlinPB10} is done with the DeepMartNet with the fully connected NN $(10, 10, 10, 1)$. The activation function for the first hidden layer is Tanh, and the one for the second hidden layer is GELU with Tanh approximation. $M_{\text{tot}} = 10^6$ paths are used and the size of mini-batch $M_b = 4000$ paths are sampled in each epoch. The total number of epoch of the training is $3000$.
The terminal time for the paths is $T=0.25$ and time step is $t=0.01$. The learning rate starts at 0.01, and decreased by a factor of 0.99 every 100 epochs.  The training takes less than 20 minutes.

\begin{figure}
    \centering
    \includegraphics[width=0.45\textwidth]{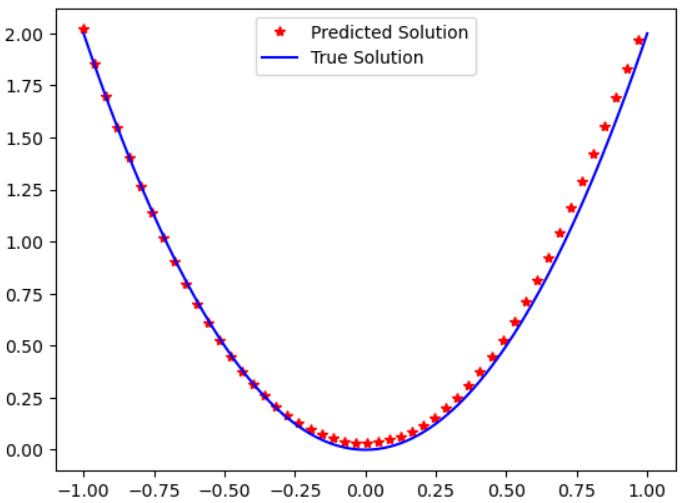}
    \includegraphics[width=0.45\textwidth]{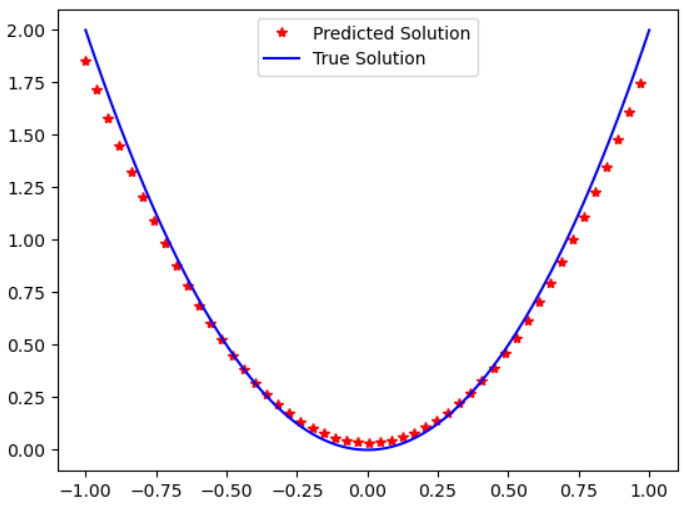}
    \includegraphics[width=0.45\textwidth]{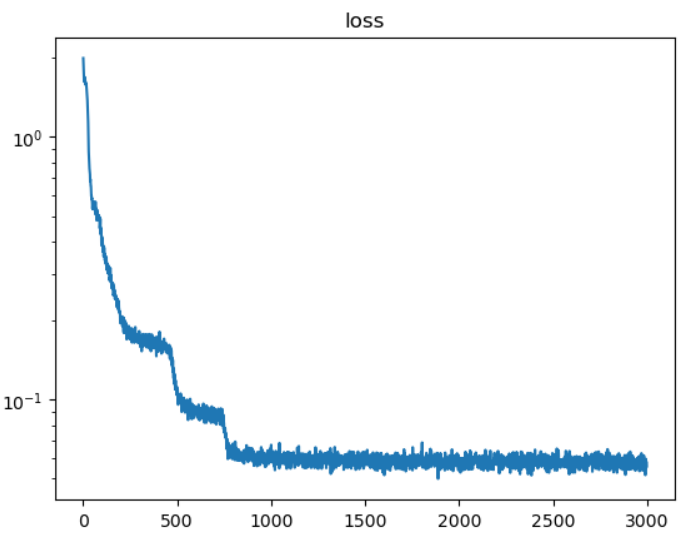}
    \includegraphics[width=0.45\textwidth]{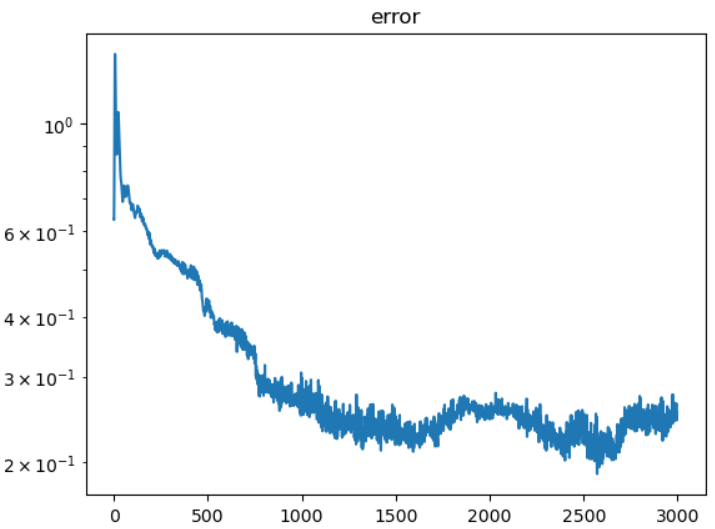}
    \caption{DeepMartNet solution of nonlinear PBE (\ref{nonlinPBsinh}) in a 10-dimensional unit ball.  Upper left: $u(x\boldsymbol{e_1})$ where $\boldsymbol{e_1} = (1, 0, \cdots, 0)$; Upper right: $u(x\boldsymbol{e})$ where $\boldsymbol{e} = (1, 1, \cdots, 1)$; Lower left: loss vs. epoch; right: relative $L^2$ error vs epoch.}
    \label{fig:NonlinPB10}
\end{figure}

\end{itemize}

\subsection{Eigenvalue problem for elliptic equations}
In this section, we will apply the DeepMartNet to solve elliptic eigenvalue problems in bounded and unbounded domains for both self-adjoint and non-self adjoint operators.
\subsubsection{Eigenvalue problem of the Laplace equation in a cube in $R^{10}$}

First, we consider a self-adjoint eigenvalue problem
\begin{align}
    -\triangle u & = \lambda u, \quad \textbf{x} \in D=[-L,L]^d, \quad d=10, \quad L=1 \\
u|_{\partial D} & = 0,
\end{align}
with the following eigenfunction
\begin{equation}
    u(x)=sin(\frac{\pi x_1}{L})\cdots sin(\frac{\pi x_d}{L}),
\end{equation}
for the lowest eigenvalue
\begin{equation}
    \lambda_1=d (\frac{\pi}{L})^2.
\end{equation}

A DeepMartNet with a fully connected structure $(10, 20, 10, 1)$ and GELU activation function are used. The total number of paths is $M = 100,000$, and for each epoch, we randomly choose $M_b = 1000$ from the paths; the time step of the paths is $\Delta t = 0.05$.  A mini-batch of $N_\text{bdry} = 2000$ points $\textbf{x}_k$ are uniformly sampled points on the boundary for every epoch of training. An Adamax optimizer with learning rate 0.05 is applied for training.
$\alpha _{\text{bdry}} = 10^3$ in \eqref{totallossEig}.  The constant $\alpha_{normal}=10, p=1, m=1, \textbf{x}_1=0, c=1$ in \eqref{lossNormal}.

To accelerate the convergence of the eigenvalue, We included an extra  term into the loss function \eqref{totallossEig}
\begin{equation}
    \alpha_{\text{eig}} |\lambda^{(i)} - \lambda^{(i - 100)}|^2,
\end{equation}
where $i$ is the index of the current epoch, which corresponds to the residual of the equation $\frac{d \lambda}{d t}=0$ as the training time $t \rightarrow \infty$, and ensures that the eigenvalue $\lambda$ will converge and the penalty constant $\alpha_{\text{eig}}$ is chosen as $2.5 \times 10 ^{-8}$ in this case. The terminal time for the paths is $T=0.6$ and time step is $t=0.01$.  The learning rate starts at 0.02, and decreased by a factor of 0.995 every 100 epochs.  The training takes less than 20 minutes.

\begin{figure}[htbp]
    \centering
    \includegraphics[width=0.45\textwidth]{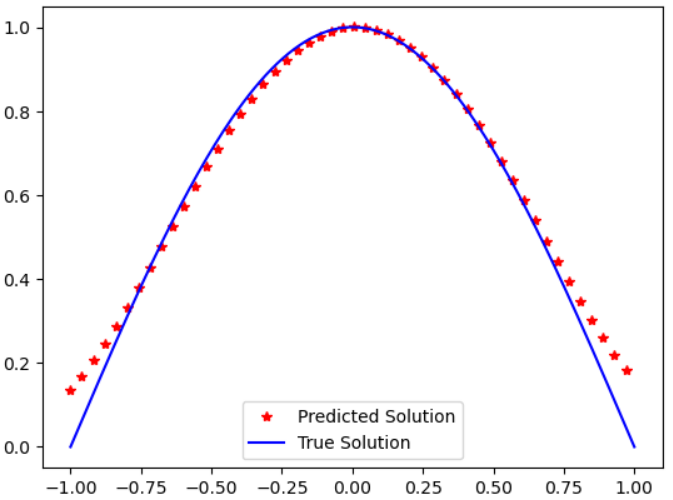}
    \includegraphics[width=0.45\textwidth]{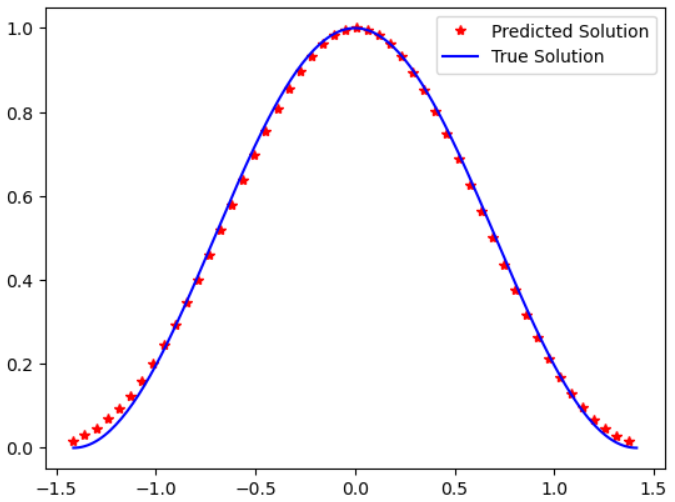}
    \includegraphics[width=0.45\textwidth]{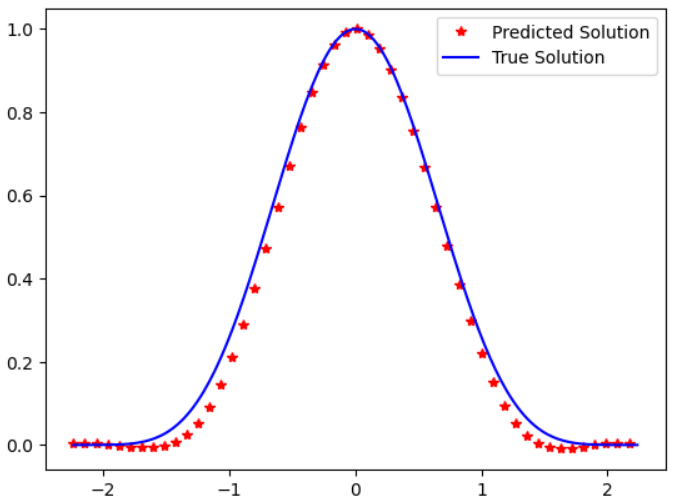}
    \includegraphics[width=0.45\textwidth]{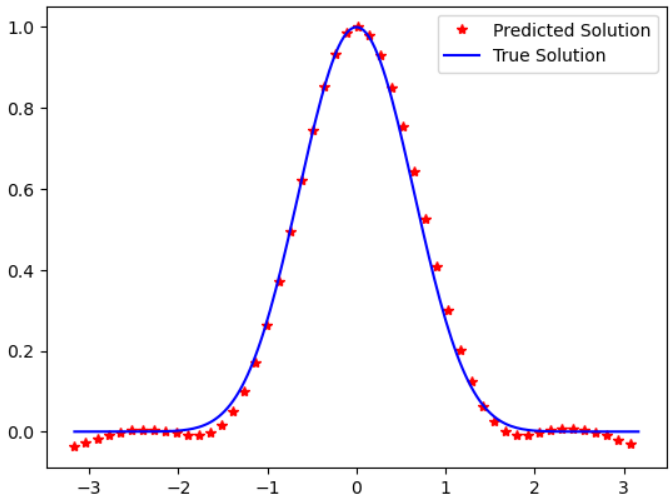}
    \includegraphics[width=0.45\textwidth]{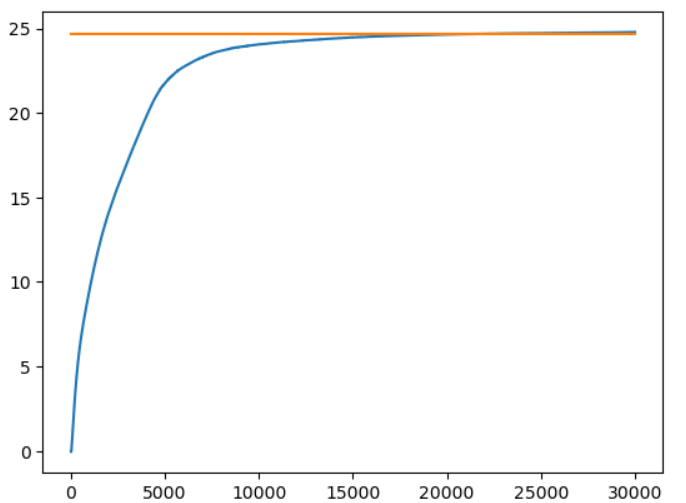}
    \includegraphics[width=0.45\textwidth]{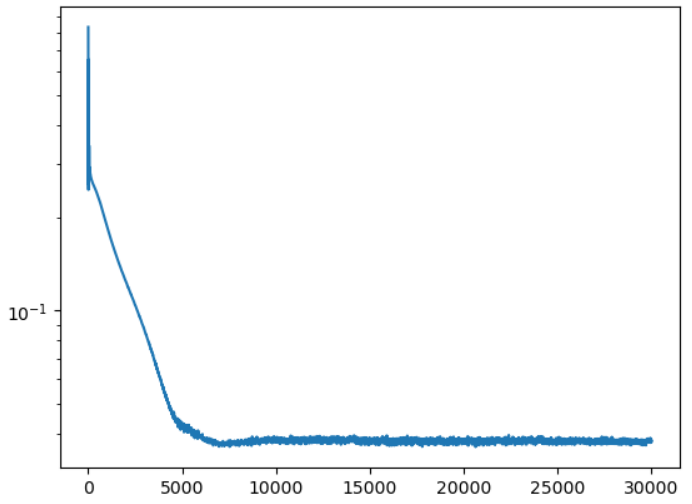}
    \caption{A numerical result for eigenvalue problem of Laplace equation in a cube $[-L,L]^{10}$ in $R^{10}$.  Upper left to middle right: true and predicted value of $u$ at $x\boldsymbol{e}$ where $\boldsymbol{e} = k^{-1/2}(\underset{k 1's}{1,\cdots,1}, 0, 0,\cdots, 0)$ for $k=1,2,5,10$
    Lower left: The predicted eigenvalue $\lambda$ versus the number of epoch, (orange) horizontal line shows the true eigenvalue; lower right: The relative error $L^2$ error versus the number of epoch. }
    \label{fig:eigLap10}
\end{figure}

\subsubsection{A non-self adjoint eigenvalue problem for the Fokker-Planck equations}
Here we consider the following a non-self adjoint eigenvalue problem of the Fokker-Planck equation
\begin{equation}
        -\Delta \psi - \nabla \cdot \left( \psi \nabla W \right) +c \psi= -\Delta \psi - \nabla W \cdot \nabla \psi - \Delta W \psi + c \psi = \lambda \psi, \quad \textbf{x} \in R^d
    \label{FokkerPlanckOriginal}
\end{equation}
where the eigenfunction for the eigenvalue $\lambda=c$ with a zero decay condition at $\infty$ is
\begin{equation}
    \psi(\textbf{x}) = e^{-W(\textbf{x})}.
\end{equation}
Here, we will consider a quadratic potential for our numerical tests $$W(\textbf{x})=||\textbf{x}||^2, \quad \textbf{x} \in R^d.$$

 Equation \eqref{FokkerPlanckOriginal} can re-written as
    \begin{equation}
    \mathcal{L} \psi = \frac{1}{2} \Delta \psi + \frac{1}{2} \nabla W \cdot \nabla \psi= -\left(  \frac{1}{2} \Delta W -\frac{1}{2}c+ \frac{1}{2}\lambda \right) \psi.
    \label{FokkerPlanck}
\end{equation}
The generator for the SDE $\mathcal{L}$ will have a drift and a diffusion as
\begin{equation*}
    \mu = \frac{1}{2} \nabla W \quad \text{and} \quad \sigma = I_{d \times d}.
    \label{FokkerPlanckMartCs}
 \end{equation*}

In order to enforce explicitly the decay condition of the eigenfunction at the infinite, a mollifier of the following form will be used as a pre-factor for the DNN solution $u_{\theta}(\textbf{x})=\rho(\textbf{x}) \tilde u_{\theta}(\textbf{x})$  to the eigenfunction,,
\begin{equation}
    \rho(\textbf{x}) = \frac{1}{1 + (\frac{||\textbf{x}||}{\alpha})^2},
    \label{mollifier}
\end{equation}
where $\alpha$ is a constant and  $\alpha = \frac{5}{11}$.

The Martingale loss of \eqref{loss_Mart_eig} for this case will be
\begin{align}
    {Loss}_{mart} = & \frac{1}{\Delta t} \frac{1}{N} \sum_{i = 0}^{N - 1}   \frac{1}{|A_i|^2} \left(\sum_{m = 1}^{|A_i|}
     u_{\theta}(\textbf{x}_{\textbf{i + 1}}^{(m)})  -
    u_{\theta}(\textbf{x}_{\textbf{i}}^{(m)}) \right. \nonumber \\
    + & (\frac{1}{2} \Delta W(\textbf{x}_{\textbf{i}}^{(m)}) -\frac{1}{2}c+ \frac{1}{2}\lambda) u_{\theta}(\textbf{x}_{\textbf{i}}^{(m)})\Delta t \Big )^2.
\end{align}


For the mini-batches in this case, we take the size of each $A_i$ to be between $\frac{M}{200}$
and $\frac{M}{25}$ as a random assortment of the $M$ total
trajectories.

To speed up the convergence of the eigenvalue and eigenfunctions, we found out that by taking fractional powers of each individual loss term is helpful, namely, the total loss is now modified as
\begin{equation}
Loss_{total}= \left( \left( Loss_{mart} \right)^p + \alpha_{normal} (Loss_{normal})^q \right)^r.
\label{epsdel}
\end{equation}
In our numerical tests, a typical choice is $p=3/8, q=1, r=3/4$, and for the normalization $Loss_{normal}$, we set $\alpha_{normal}=50$ and $c=30$, $m=2, \textbf{x}_1=0, \textbf{x}_2=(1, 0, \cdots, 0)$ in \eqref{lossNormal}.

\begin{figure}[htbp]
  \centering
    \includegraphics[scale = 0.29]{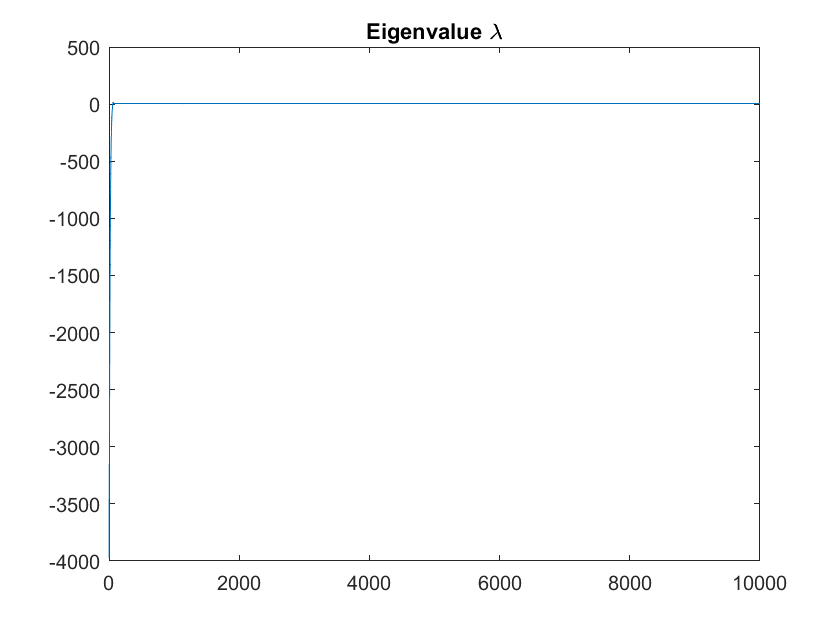}
    \includegraphics[scale = 0.29]{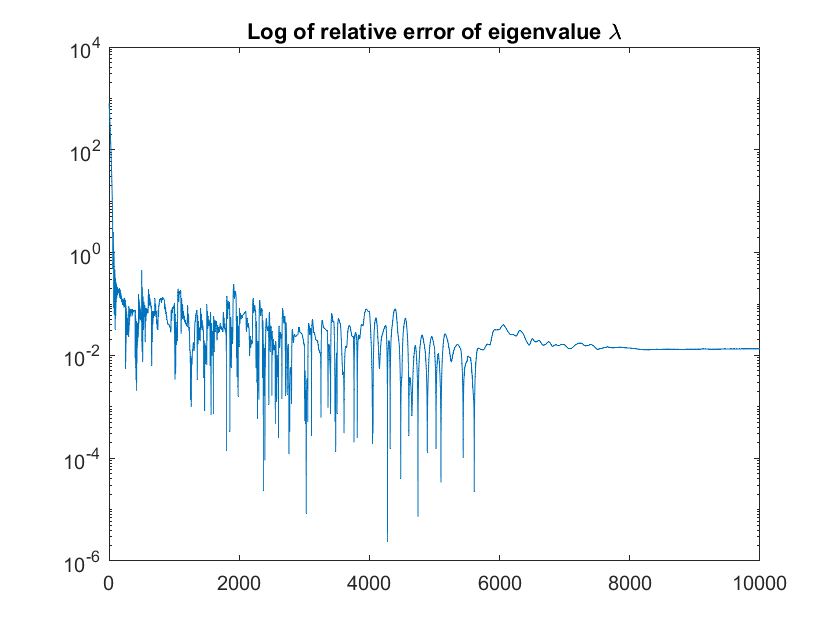}
    \includegraphics[scale = 0.29]{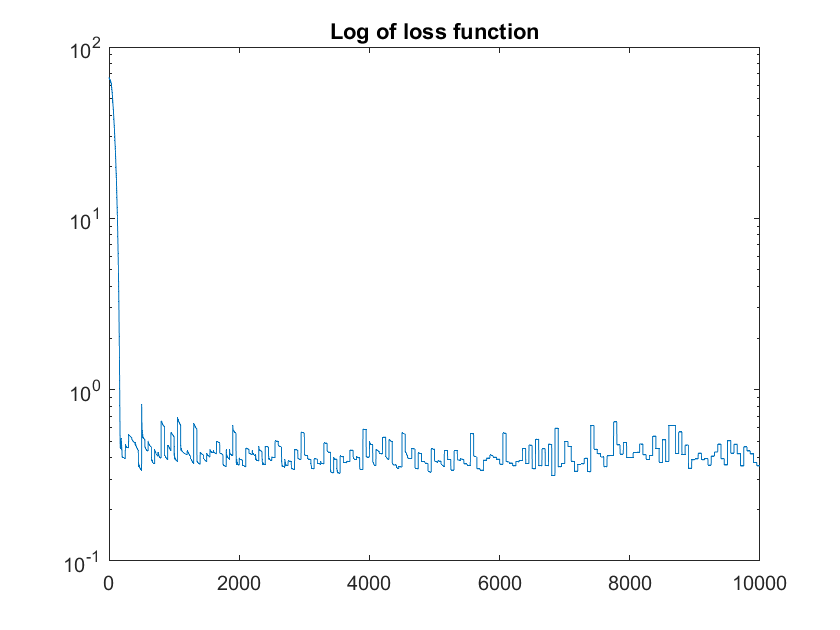}
    \includegraphics[scale = 0.29]{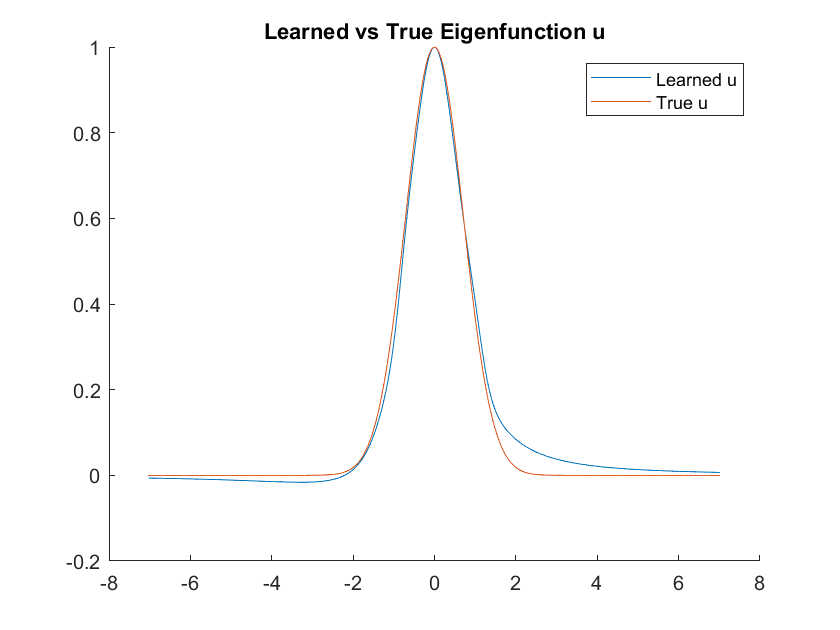}

    \caption{ Eigen-value problem of Fokker-Planck equation in $R^d$  $d=5$ for eigenvalue $\lambda = 5$.
    (Top left) Convergence of eigenvalue, (top right) history of eigenvalue error, (bottom left) history of loss function, (bottom right) Learn and exact eigenfunction along the diagonal of the domain.}
    \label{fpev_d5}
\end{figure}

\begin{figure}[htbp]
    \centering
    \includegraphics[scale = 0.29]{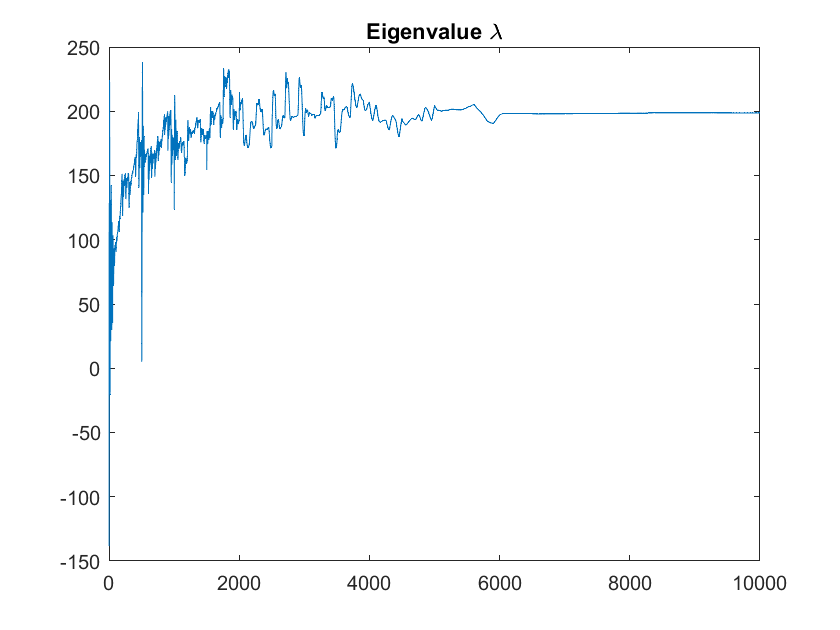}
        \includegraphics[scale = 0.29]{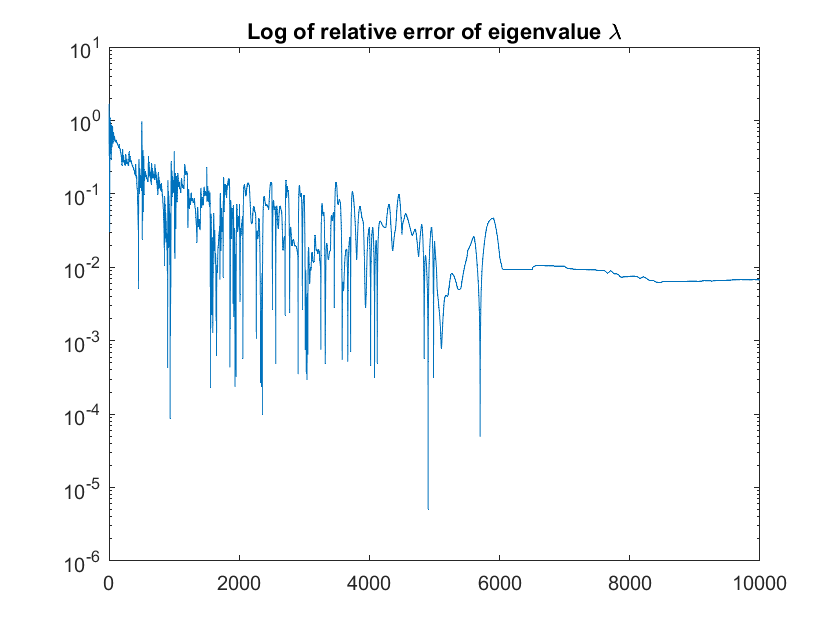}
    \includegraphics[scale = 0.29]{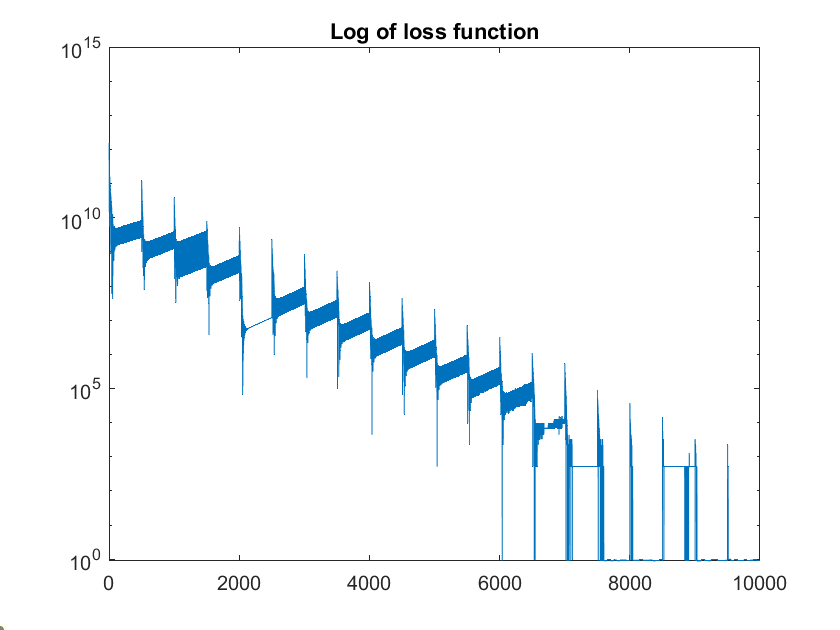}
    \includegraphics[scale = 0.29]{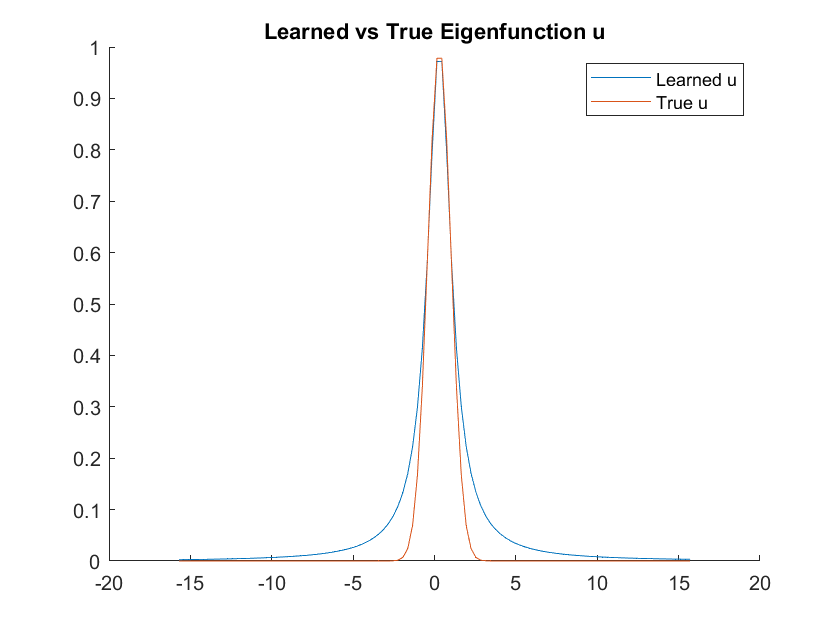}


\caption{ Eigenvalue problem of Fokker-Planck equation in $R^d$  $d=200$ for eigenvalue $\lambda = 200$.
(Top left) Convergence of eigenvalue, (top right) history of eigenvalue error, (bottom left) history of loss function, (bottom right) Learn and exact eigenfunction along the diagonal of the domain.}

    \label{fpev_d200}
\end{figure}

Figs. \ref{fpev_d5}  and \ref{fpev_d200} show the learned eigenvalues ($\lambda=5$ for $d=5$
and $\lambda=200$ for $d=200$) and eigenfunctions, respectively. The  $k=3$-trapezoidal rule is used in the Martingale loss ${Loss}_{mart}$ and the other numerical parameters used are listed as follows
\begin{itemize}
    \item The total number of paths starting from the origin $M=9,000,   24,000$ for $d=5, 200$, respectively.
    \item The number of time steps $N=1350, 1300$ and the terminal time $T=9$ for $d=5, 200$, respectively.
    \item the learning rate is 1/150 and is halved every 500 epochs starting at epoch 500 and
halved twice at epoch 7500 for stabilization.
\item A  In the following numerical tests, a fully connected network with layers (d, 6d, 3d, 1) with a Tanh activation function is used for the eigenfunction while a fully connected network (1, d, 1) with a $ReLu^9$ activation function with a constant value input  is used to represent the eigenvalue.
\item An Adamax optimizer is applied for training.
\end{itemize}

The final relative error in the eigenvalues after 10,000 epochs of training is $1.3 \times 10^{-2},  6.7 \times 10^{-3} $ for $d=5, 200$, respectively. And the relative $L^2$ error of the eigenfunction calculated along the diagonal of the domain is $2.6 \times 10^{-2},   2.9 \times 10^{-2} $ for $d=5, 200$, respectively. The training takes 25 minutes for the case of $d=200$.

\section{Conclusion}

In this paper, we introduced a Martingale based neural network, DeepMartNet, for solving the BVPs and eigenvalue problems of elliptic operators with Dirichlet boundary conditions. The DeepMartNet enforces the Martingale property for the PDE solutions through the stochastic gradient descent (SGD) optimization of the Martingale property loss functions. The connection between the Martingale definition and the mini-batches used in stochastic gradient computation shows a natural fit between the SGD optimization in DNN training and the Martingale problem formulation of the PDE solutions.
Numerical results in high dimensions for both BVPs and eigenvalue problems show the promising potential in approximating high dimensional PDE solutions. The numerical results show that the DeepMartNet  extracts more information for the PDEs solution over the whole solution domain than the traditional one point solution Feynman-Kac formula from the same set of diffusion paths originating from just one point.

Future work on DeepMartNet will include PDEs with Neumann and Robin boundary conditions where the underlying diffusion will be a reflecting one and the local time of the reflecting diffusion will be computed and included in the definition of the Martingale loss function. Another area of application is in optimal stochastic control where the Martingale optimality condition for the control and the backward SDE for the value function can be used to construct the loss function for the control as well as the value function \cite{deepMart}. Another important area for the application of the DeepMartNet is to solve low-dimension PDEs in complex geometries as the method uses diffusion paths to explore the domain and therefore can handle highly complex geometries such as the interconnects in microchip designs and nano-particles in material sciences and molecules in biology. Finally, The convergence analysis of the DeepMartNet, especially to a given eigenvalue,  and the choices of mini-batches of diffusion paths and various hyper-parameters, which can affect the convergence of the learning strongly, in the loss functions are important issues to be addressed.

 \section*{Acknowledgement} W. C. would like to thank Elton Hsu and V. Papanicolaou for the helpful discussion about their work on probabilistic solutions of Neumann and Robin problems.


\end{document}